\documentclass[acmsmall]{acmart}
\acmJournal{TOIS}

\usepackage{mathtools,bm}
\usepackage{booktabs,tabularx,longtable,array}
\usepackage{pdflscape}
\usepackage{enumitem}
\usepackage{ragged2e}
\usepackage{forest}
\usepackage{adjustbox}

\newcolumntype{Y}{>{\RaggedRight\arraybackslash}X}
\newcolumntype{P}[1]{>{\RaggedRight\arraybackslash}p{#1}}
\newcommand{\R}{\mathbb{R}}
\newcommand{\B}{\{0,1\}}
\newcommand{\tr}{\operatorname{tr}}

\newcommand{\AppField}[1]{\par\smallskip\noindent\textbf{#1}}
\acmYear{2026}
\copyrightyear{2026}
\acmMonth{8}
\acmVolume{0}
\acmNumber{0}
\acmArticle{0}
\acmPrice{15.00}
\acmDOI{10.0000/0000000}

\title{Mathematical Programming in Machine Learning and Artificial Intelligence: A Unified Taxonomy of Models and Applications}

\author{Chaosheng Dong}
  \email{ensteindcs@gmail.com}

\begin{abstract}
Mathematical programming provides a common language for many decisions embedded in modern machine-learning (ML) and artificial-intelligence (AI) systems: selecting retrieval context, routing tokens, allocating inference compute, fitting structured predictors, protecting against distribution shift, and balancing competing objectives.  However, the relevant literature is fragmented across optimization, information retrieval, recommendation, natural-language processing, computer vision, and learning theory.  This paper organizes  various applications under common mathematical programming paradigms: linear, quadratic, binary and mixed-integer, conic, bilevel, multi-objective, inverse, distributionally robust, submodular, and min--max optimization.  We normalize the models with a mostly unified notation and, for every application, identify inputs, decision variables, a principal formulation, structural properties, solution strategies, and limitations.  Across paradigms, we compare tractability, relaxation quality, decomposition, approximation guarantees, and scalability bottlenecks.  The paper shows that mathematical programming is most useful not as a claim that all learning is LP or MIP, but as a disciplined interface between predictions and constrained decisions.
\end{abstract}


\begin{document}
\maketitle

\section{Introduction}
The output of an ML model is often  an input to a decision.  A retrieval score must be converted into a context that fits a token window; a mixture-of-experts gate must respect capacity; a recommender must construct a list rather than score items independently; and an alignment procedure must remain useful under shifts in preference data.  These decisions couple predictions through budgets, assignments, diversity, dynamics, or adversaries.  Mathematical programming (MP) makes those couplings explicit and exposes which guarantees survive when learned quantities enter an optimization model.

This perspective is especially natural for information systems.  Search, recommendation, retrieval-augmented generation (RAG), notification delivery, and large language model (LLM) inference all combine statistical estimates with resource-constrained actions.  The same modeling vocabulary also clarifies classical learning procedures, including sparse recovery, manifold learning, grouped regression, structured model selection, and adversarial generation.  The point is not that every algorithm is implemented by a generic solver.  Rather, an MP formulation can serve as an exact computational model, an equivalent representation that reveals structure, or a normative reformulation useful for analysis and system design.

This paper makes four contributions.  First, it provides a taxonomy connecting MP paradigms to ML/AI applications.  Second, it gives normalized formulations using a mostly unified notation, so that similarities such as budgets, simplex constraints, and nested objectives are visible across domains.  Third, it distinguishes \emph{explicit}, \emph{equivalent}, and \emph{natural} MP relationships to the cited primary literature.  Fourth, it compares tractability, relaxations, decomposition, approximation, and scalability, including where solver certificates are practical and where a formulation primarily guides a heuristic.

The paper is organized as follows. Sections~\ref{sec:scope}--\ref{sec:minmax} develop the paradigms and their application models; Section~\ref{sec:comparison} compares them; and Section~\ref{sec:challenges} identifies research directions.

%

\begingroup

\forestset{
  mpml taxonomy/.style={
    for tree={
      grow'=0,
      anchor=west,
      parent anchor=east,
      child anchor=west,
      align=left,
      font=\scriptsize,
      inner xsep=1.0pt,
      inner ysep=0.25pt,
      outer sep=0pt,
      l sep=5.5mm,
      s sep=0.30mm,
      edge={draw=black,line width=0.35pt},
      edge path'={
        (!u.parent anchor) -- ++(0.9mm,0) |- (.child anchor)
      },
    },
  },
  mpml root/.style={
    draw=black,
    line width=0.45pt,
    text width=2.70cm,
    font=\bfseries\scriptsize,
    align=center,
    inner xsep=2.5pt,
    inner ysep=2.0pt,
  },
  mpml family/.style={
    text width=2.55cm,
    font=\itshape\scriptsize,
    for children={s sep=1.0mm},
  },
  mpml paradigm/.style={
    text width=3.35cm,
    font=\bfseries\scriptsize,
    for children={s sep=0.30mm},
  },
  mpml subclass/.style={
    text width=2.40cm,
    font=\itshape\scriptsize,
  },
  mpml application/.style={
    text width=6.00cm,
    font=\normalfont\scriptsize,
  },
  mpml deep application/.style={
    text width=5.00cm,
    font=\normalfont\scriptsize,
  },
}

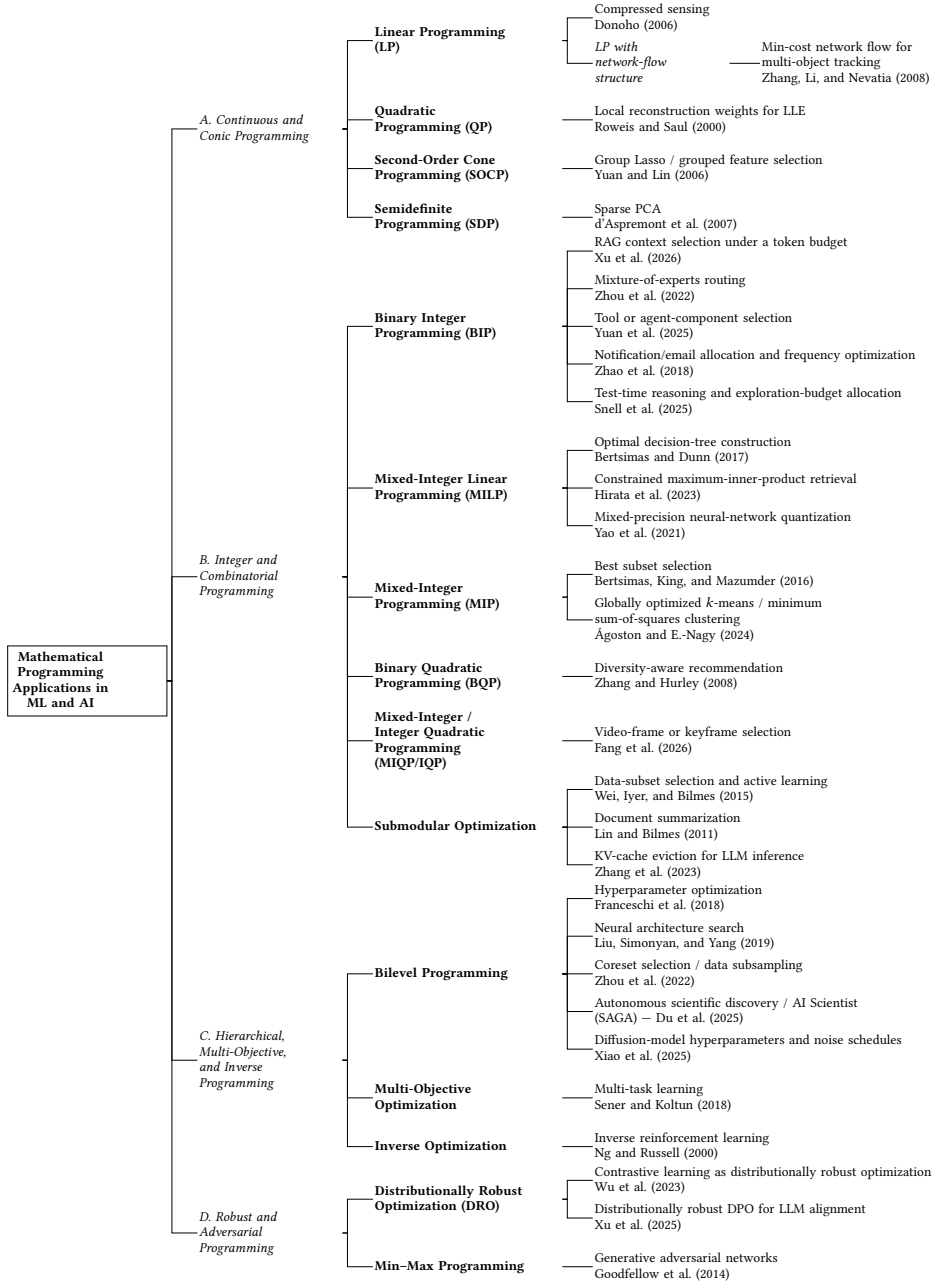
\begin{figure*}[t]
  \centering
  \begin{adjustbox}{max totalsize={\textwidth}{0.84\textheight},center}
  \begin{forest}
    mpml taxonomy
    [Mathematical\\Programming\\Applications in\\ML and AI,
      mpml root,
      for children={s sep=3.0mm}
      [A. Continuous and\\Conic Programming, mpml family
        [Linear Programming\\(LP), mpml paradigm
          [Compressed sensing\\
            Donoho (2006), mpml application]
          [LP with\\network-flow\\structure, mpml subclass
            [{Min-cost network flow for\\multi-object tracking\\
              Zhang, Li, and Nevatia (2008)}, mpml deep application]
          ]
        ]
        [Quadratic\\Programming (QP), mpml paradigm
          [Local reconstruction weights for LLE\\
           Roweis and Saul (2000), mpml application]
        ]
        [Second-Order Cone\\Programming (SOCP), mpml paradigm
          [Group Lasso / grouped feature selection\\
           Yuan and Lin (2006), mpml application]
        ]
        [Semidefinite\\Programming (SDP), mpml paradigm
          [Sparse PCA \\
          d'Aspremont et al. (2007), mpml application]
        ]
      ]
      [B. Integer and\\Combinatorial\\Programming, mpml family
        [Binary Integer\\Programming (BIP), mpml paradigm
          [RAG context selection under a token budget\\
           Xu et al. (2026), mpml application]
          [Mixture-of-experts routing\\
           Zhou et al. (2022), mpml application]
          [Tool or agent-component selection\\
           Yuan et al. (2025), mpml application]
          [Notification/email allocation and frequency optimization \\
            Zhao et al. (2018), mpml application]
          [Test-time reasoning and exploration-budget allocation \\
            Snell et al. (2025), mpml application]
        ]
        [Mixed-Integer Linear\\Programming (MILP), mpml paradigm
          [Optimal decision-tree construction\\
           Bertsimas and Dunn (2017), mpml application]
          [Constrained maximum-inner-product retrieval\\
           Hirata et al. (2023), mpml application]
          [Mixed-precision neural-network quantization\\
           Yao et al. (2021), mpml application]
        ]
        [Mixed-Integer\\Programming (MIP), mpml paradigm
          [{Best subset selection\\
            Bertsimas, King, and Mazumder (2016)}, mpml application]
          [Globally optimized $k$-means / minimum\\
           sum-of-squares clustering\\
           \'{A}goston and E.-Nagy (2024),
           mpml application]
        ]
        [Binary Quadratic\\Programming (BQP), mpml paradigm
          [Diversity-aware recommendation\\
           Zhang and Hurley (2008), mpml application]
        ]
        [Mixed-Integer /\\Integer Quadratic\\Programming\\
         (MIQP/IQP), mpml paradigm
          [Video-frame or keyframe selection\\
           Fang et al. (2026), mpml application]
        ]
        [Submodular Optimization, mpml paradigm
          [{Data-subset selection and active learning\\
            Wei, Iyer, and Bilmes (2015)}, mpml application]
          [Document summarization\\
           Lin and Bilmes (2011), mpml application]
          [KV-cache eviction for LLM inference\\
           Zhang et al. (2023), mpml application]
        ]
      ]
      [{C. Hierarchical,\\Multi-Objective,\\and Inverse\\
        Programming}, mpml family
        [Bilevel Programming, mpml paradigm
          [Hyperparameter optimization\\
           Franceschi et al. (2018), mpml application]
          [{Neural architecture search\\
            Liu, Simonyan, and Yang (2019)}, mpml application]
          [Coreset selection / data subsampling\\
           Zhou et al. (2022), mpml application]
          [Autonomous scientific discovery / AI Scientist\\
           (SAGA) --- Du et al. (2025), mpml application]
          [Diffusion-model hyperparameters and noise schedules \\
            Xiao et al. (2025), mpml application]
        ]
        [Multi-Objective\\Optimization, mpml paradigm
          [Multi-task learning\\
           Sener and Koltun (2018), mpml application]
        ]
        [Inverse Optimization, mpml paradigm
          [Inverse reinforcement learning\\
           Ng and Russell (2000), mpml application]
        ]
      ]
      [D. Robust and\\Adversarial\\Programming, mpml family
        [Distributionally Robust\\Optimization (DRO), mpml paradigm
          [Contrastive learning as distributionally robust optimization \\
          Wu et al. (2023), mpml application]
          [Distributionally robust DPO for LLM alignment\\
           Xu et al. (2025), mpml application]
        ]
        [Min--Max Programming, mpml paradigm
          [Generative adversarial networks\\
           Goodfellow et al. (2014), mpml application]
        ]
      ]
    ]
  \end{forest}
  \end{adjustbox}
  \caption{Taxonomy of mathematical programming applications in machine
    learning and artificial intelligence.}
  \Description{A left-to-right hierarchical line-based taxonomy with four
    organizational families, mathematical-programming paradigms, and thirty
    machine-learning and artificial-intelligence applications, each paired
    with a representative author--year citation. Network-flow optimization is
    shown as a structured subclass of linear programming.}
  \label{fig:mp-ml-taxonomy}
\end{figure*}

\endgroup

\section{Scope, Taxonomy, and Mathematical Programming Preliminaries}
\label{sec:scope}
In operations research, MP denotes finite-dimensional optimization with  objectives and constraints, especially linear programming (LP), quadratic programming (QP), mixed-integer programming (MIP), and conic programming.  Our  scope in this paper is deliberately broader: it also includes bilevel, multi-objective, inverse, distributionally robust optimization (DRO), submodular optimization, and min-max games.  These classes retain the central MP discipline of explicit decisions, feasible sets, and objectives, even when they are nonsmooth, combinatorial, infinite-dimensional before dualization, or game-theoretic.  This broader scope is necessary to describe modern ML/AI faithfully; it does not imply that every class reduces efficiently to LP or MIP.

We use three relationship labels.  \emph{Explicit formulation} means that the cited work directly formulates the application in the stated paradigm or an unmistakably matching form.  \emph{Equivalent formulation} means that the displayed program is mathematically equivalent to the cited objective but may expose auxiliary variables or a lifted representation not emphasized by the source.  \emph{Natural MP reformulation} means that we translate a source's algorithmic or system description into a useful MP model; we do not make the claim that the paper solved that displayed program.

For a generic continuous model, let the decision be $\mathbf{x}\in\mathcal X\subseteq\R^d$, the objective be $f(\mathbf{x})$, inequality functions be $g_u(\mathbf{x})$ for $u$ in a finite index set $\mathcal U$, and equality functions be $h_v(\mathbf{x})$ for $v$ in a finite index set $\mathcal W$.  These symbols are local to the next display:
\begin{equation}
  \min_{\mathbf{x}\in\mathcal X} f(\mathbf{x})
  \quad\text{s.t.}\quad
  g_u(\mathbf{x})\le 0\ (u\in\mathcal U),\qquad
  h_v(\mathbf{x})=0\ (v\in\mathcal W).
  \label{eq:generic-mp}
\end{equation}
Convexity of $f$, the functions $g_u$, and the set $\mathcal X$, together with linear equalities, makes local optima global.  Binary variables replace a continuous feasible region by a discrete one, usually making worst-case solution time exponential unless special structure applies.  Relaxations enlarge a feasible set to produce bounds; rounding or primal heuristics recover implementable decisions.  Decomposition separates weakly coupled blocks, while approximation algorithms trade exactness for a provable or empirical quality guarantee.

\subsection{Taxonomy at a Glance}
We report the taxonomy of mathematical programming applications in machine
learning and artificial intelligence in Table~\ref{tab:taxonomy}.

\begin{landscape}
\scriptsize
\setlength{\tabcolsep}{3pt}
\begin{longtable}{P{0.4cm}P{1.7cm}P{3.5cm}P{2.7cm}P{4.3cm}P{2.4cm}P{1.9cm}}
\caption{Unified taxonomy. }\label{tab:taxonomy}\\
\toprule
No. & Paradigm & ML/AI application & Main decision variable & Principal constraint or structure & Representative citation & MP relationship\\
\midrule
\endfirsthead
\multicolumn{7}{l}{\tablename~\thetable\ continued}\\
\toprule
No. & Paradigm & ML/AI application & Main decision variable & Principal constraint or structure & Representative citation & MP relationship\\
\midrule
\endhead
\midrule
\multicolumn{7}{r}{Continued on next page}\\
\endfoot
\bottomrule
\endlastfoot
1 & LP & Compressed sensing & Sparse signal $\mathbf{x}$, absolute-value epigraph & Linear measurements; $\ell_1$ epigraph & \citet{donoho2006compressed} & Equivalent formulation\\
2 & LP (network flow) & Min-cost network flow for multi-object tracking & Detection and transition flows & Flow conservation, capacities, source/sink & \citet{zhang2008tracking} & Explicit formulation\\
3 & QP & Local reconstruction weights in LLE & Neighbor weights $w_{ij}$ & linear sum-to-one; neighborhood support & \citet{roweis2000lle} & Explicit formulation\\
4 & BIP/MMKP & RAG context selection & Passage indicators $z_i$ & Token/resource budgets; group choice & \citet{xu2026selfcorrecting} & Explicit formulation\\
5 & BIP & Mixture-of-experts routing & Token--expert assignments $z_{ik}$ & Per-token assignment; expert capacity & \citet{zhou2022expertchoice} & Natural MP reformulation\\
6 & BIP & Tool or agent-component selection & Component indicators $z_i$ & Budget and incompatibility constraints & \citet{yuan2025agentcomposition} & Explicit formulation\\
7 & BIP & Notification/email allocation and frequency & User--message--time indicators & Frequency, channel, quota, and global budgets & \citet{zhao2018notification,zhang2020billion} & Natural MP reformulation\\
8 & BIP/MCKP & Test-time reasoning and exploration budget & Query--compute-level choices & Exactly one level; global compute & \citet{snell2025testtime} & Natural MP reformulation\\
9 & MILP & Optimal decision-tree construction & Splits, leaf assignments, predictions & Tree logic and misclassification linearization & \citet{bertsimas2017oct} & Explicit formulation\\
10 & MIQP & Best subset selection & Coefficients $\theta_j$, indicators $z_j$ & Cardinality and coefficient linking & \citet{bertsimas2016subset} & Explicit formulation\\
11 & MILP/MIQP & Constrained maximum-inner-product retrieval & Retrieved-item indicators $z_i$ & Cardinality, coverage/fairness; pairwise diversity & \citet{hirata2023categorical} & Natural MP reformulation\\
12 & BIP/ILP & Mixed-precision neural-network quantization & Layer--bit-width assignments & Exactly one width; size or latency budget & \citet{yao2021hawq} & Explicit formulation\\
13 & BQP/MIQP & Diversity-aware recommendation & List indicators $z_i$ & Cardinality; pairwise similarity penalty & \citet{zhang2008monotony} & Natural MP reformulation\\
14 & MIQP/\newline MINLP/\newline MILP & Globally optimized $k$-means clustering & Assignments $z_{ik}$, centroids $\bm\mu_k$ & Exactly one cluster; assignment--centroid coupling & \citet{agoston2024kmeans,burgard2023clustering} & Explicit formulation\\
15 & BQP/IQP & Video frame or keyframe selection & Frame indicators $z_i$ & Cardinality/temporal budget; pairwise redundancy & \citet{fang2026keyframe} & Explicit formulation\\
16 & SOCP & Group Lasso & Coefficients $\bm\beta_g$, epigraphs $t_g$ & Group-norm cones & \citet{yuan2006grouplasso} & Equivalent formulation\\
17 & SDP & Sparse PCA & Lifted matrix $\mathbf{X}\succeq0$ & Trace normalization; rank relaxation; sparsity & \citet{daspremont2004sparsepca} & Explicit formulation\\
18 & Bilevel & Hyperparameter optimization & Hyperparameters $\lambda$, weights $\theta$ & Validation outer/training inner & \citet{franceschi2018bilevel} & Explicit formulation\\
19 & Bilevel & Neural architecture search & Architecture $\lambda$, weights $\theta$ & Relaxed architecture outer/weight-training inner & \citet{liu2019darts} & Explicit formulation\\
20 & Bilevel & Coreset selection and data subsampling & Selection probabilities/indicators and weights & Coreset-size outer; selected-data training inner & \citet{zhou2022coreset} & Explicit formulation\\
21 & Bilevel & Autonomous scientific discovery & Objective design $\lambda$, candidate design $\mathbf{x}$ & Goal-evolution outer; scientific search inner & \citet{du2025saga,ma2024scientific} & Equivalent formulation\\
22 & Bilevel & Diffusion hyperparameter/noise-schedule optimization & Schedule $\lambda$, diffusion parameters $\theta$ & Quality outer; training/inference inner & \citet{xiao2025diffusion} & Explicit formulation\\
23 & Multi-objective & Multi-task learning & Shared/task model parameters $\theta$ & Pareto order; simplex gradient QP & \citet{sener2018multitask} & Explicit formulation\\
24 & Inverse optimization & Inverse reinforcement learning & Reward weights $\theta$ & Expert policy optimality/Bellman inequalities & \citet{ng2000irl} & Explicit formulation\\
25 & DRO & Contrastive learning as DRO & Encoder parameters $\theta$, negative law $Q$ & Divergence ambiguity around empirical negatives & \citet{wu2023contrastive} & Explicit formulation\\
26 & DRO & Distributionally robust DPO & Policy parameters $\theta$, preference law $Q$ & Wasserstein or KL ambiguity set & \citet{xu2025robustdpo} & Explicit formulation\\
27 & Submodular & Data subset selection and active learning & Selected set $S$ & Cardinality; diminishing returns & \citet{wei2015submodularity} & Explicit formulation\\
28 & Submodular & Document summarization & Selected sentence set $S$ & Word budget; coverage and diversity & \citet{lin2011summarization} & Explicit formulation\\
29 & Submodular & KV-cache eviction for LLM inference & Retained token set $S_\tau$ & Dynamic cache budget; heavy hitters plus recency & \citet{zhang2023h2o} & Explicit formulation\\
30 & Min--max & Generative adversarial networks & Generator $G$, discriminator $D$ & Two-player saddle objective & \citet{goodfellow2014gan} & Explicit formulation\\
\end{longtable}
\end{landscape}

\subsection{Unified Notation}
\label{sec:notation}
Bold lowercase letters denote vectors, bold uppercase letters matrices, calligraphic letters sets, and ordinary lowercase letters scalars.  A symbol may acquire a narrowly specialized meaning only when that meaning is defined next to its equation.  Table~\ref{tab:notation} lists the  notations.

\begin{table*}
\caption{List of notation. }
\label{tab:notation}
\scriptsize
\begin{tabularx}{\textwidth}{P{2.1cm}P{3.2cm}YP{3.0cm}}
\toprule
Symbol & Domain or type & Meaning & First section used\\
\midrule
$i,j$ & indices & Objects, candidates, samples, or items & Linear Programming\\
$q$ & index & Query or request & Binary and Mixed-Integer Programming\\
$g$ & index & Group & Binary and Mixed-Integer Programming\\
$k$ & index & Cluster, expert, component, or discrete option & Binary and Mixed-Integer Programming\\
$\tau$ & nonnegative integer & Time step & Linear Programming\\
$d,n,K,T$ & positive integers & Feature dimension, sample count, cluster/option count, and task/time horizon & Linear Programming\\
$\mathcal I,\mathcal J,\mathcal G$ & finite sets & Corresponding object, candidate, and group index sets & Linear Programming\\
$\mathcal K,\mathcal Q,\mathcal T$ & finite sets & Options/experts, queries, and tasks/time steps & Binary and Mixed-Integer Programming\\
$\mathcal V$ & finite set & Ground set for subset selection & Submodular Optimization\\
$\mathbf a_i$ & $\R^d$ & Observed input or feature vector & Quadratic Programming\\
$y_i$ & label/response space & Observed label or response & Binary and Mixed-Integer Programming\\
$\mathbf A,\mathbf b$ & matrix/vector & Linear operator and observed right-hand side & Linear Programming\\
$\theta$ & parameter space $\Theta$ & Learned model parameters & Binary and Mixed-Integer Programming\\
$\lambda$ & scalar or vector & Hyperparameters or trade-off parameters & Binary and Mixed-Integer Programming\\
$\mathbf x$ & $\R^d$ & Continuous decision or latent state & Linear Programming\\
$z_i,z_{ik}$ & $\{0,1\}$ & Selection or assignment decisions & Binary and Mixed-Integer Programming\\
$w_i,w_{ij}$ & $\R$ or $\R_+$ & Continuous weights & Quadratic Programming\\
$r_i$ & $\R$ & Relevance, reward, or utility & Binary and Mixed-Integer Programming\\
$c_i$ & $\R_+$ & Resource cost & Binary and Mixed-Integer Programming\\
$s_{ij}$ & $\R$ & Similarity or pairwise interaction & Binary and Mixed-Integer Programming\\
$B$ & $\R_+$ & Resource budget & Binary and Mixed-Integer Programming\\
$m$ & positive integer & Selection-cardinality limit & Binary and Mixed-Integer Programming\\
$S\subseteq\mathcal V$ & set decision & Selected subset of a ground set & Submodular Optimization\\
$L_{\mathrm{tr}}$ & real-valued function & Training loss & Bilevel Programming\\
$L_{\mathrm{val}}$ & real-valued function & Validation loss & Bilevel Programming\\
$\widehat P$ & probability measure & Empirical distribution & Distributionally Robust Optimization\\
$Q$ & probability measure & Adversarial or alternative distribution & Distributionally Robust Optimization\\
$\mathcal U(\widehat P)$ & set of measures & Ambiguity set around $\widehat P$ & Distributionally Robust Optimization\\
$\Phi,\Psi$ & real-valued functions & Upper- and lower-level objectives & Bilevel Programming\\
$F(S)$ & $\R$ & Submodular set function & Submodular Optimization\\
$\ell(\theta;\xi)$ & $\R$ & Loss on random observation $\xi$ & Distributionally Robust Optimization\\
$\mathbf 1,\mathbf I$ & vector/matrix & All-ones vector and identity matrix of contextual size & Linear Programming\\
$\lVert\cdot\rVert_p$ & norm & $\ell_p$ vector norm & Linear Programming\\
$\langle\cdot,\cdot\rangle$ & inner product & Euclidean/Frobenius inner product & Binary and Mixed-Integer Programming\\
$\mathbf X\succeq0$ & symmetric matrix & Positive-semidefinite constraint & Conic Programming\\
$\bm\epsilon$ & random vector & Exogenous noise & Bilevel Programming\\
\bottomrule
\end{tabularx}
\end{table*}

\section{Linear Programming}
\label{sec:lp}
LP optimizes an  objective over linear equalities and inequalities.  Its global optimality certificates and mature decomposition methods make it attractive when a learning-system decision is continuous or when an integral combinatorial model has an exact LP relaxation.

\subsection{Compressed Sensing}
\AppField{Context.}
Compressed sensing reconstructs a high-dimensional signal from fewer linear measurements than ambient dimension by exploiting sparsity.  Basis pursuit chooses, among signals consistent with the measurements, one with minimum $\ell_1$ norm.

\AppField{Inputs and decision variables.}
Let $\mathbf A\in\R^{p\times d}$ be a known sensing matrix with $p<d$, let $\mathbf b\in\R^p$ be the measured vector, and let the continuous decision $\mathbf x\in\R^d$ be the reconstructed signal.  With auxiliary decisions $\mathbf t\in\R_+^d$, the canonical and LP forms are
\begin{align}
  \min_{\mathbf x\in\R^d}\quad &\lVert\mathbf x\rVert_1
  &&\text{s.t.}\quad \mathbf A\mathbf x=\mathbf b,
  \label{eq:basis-pursuit}\\
  \min_{\mathbf x\in\R^d,\,\mathbf t\in\R_+^d}\quad
  &\sum_{j=1}^{d}t_j
  &&\text{s.t.}\quad
  \mathbf A\mathbf x=\mathbf b,\quad -t_j\le x_j\le t_j\quad(j=1,\ldots,d).
  \label{eq:basis-pursuit-lp}
\end{align}

\AppField{Interpretation.}
Measurement consistency is enforced exactly, while the epigraph inequalities make $t_j\ge |x_j|$; minimization makes them tight.  Thus \eqref{eq:basis-pursuit-lp} is equivalent to \eqref{eq:basis-pursuit}.  Sparsity is encouraged because the $\ell_1$ unit ball has coordinate-aligned extreme points.

\AppField{MP classification and structure.}
Equation~\eqref{eq:basis-pursuit-lp} is an LP and hence convex.  Recovery of the sparsest feasible signal is exact only under conditions on $\mathbf A$; LP optimality alone does not establish statistical identifiability.

\AppField{Solution approaches.}
Interior-point or simplex methods apply directly; first-order primal--dual, operator-splitting, and homotopy methods exploit large sparse sensing operators.  

\AppField{MP relationship and citation.} \textbf{Equivalent formulation.}  \citet{donoho2006compressed} establishes the compressed-sensing role of $\ell_1$ recovery.  Our auxiliary-variable LP is exact; we do not attribute the elementary epigraph conversion as the paper's contribution.

\AppField{Limitations.}
The model abstracts measurement noise, structured sparsity, model mismatch, and the computational cost of applying $\mathbf A$; it also treats sparsity as the relevant inductive bias without validating it for a particular signal source.

\subsection{LP with Network-Flow Structure}
\label{sec:lp-flow}
Network-flow optimization is a structured subclass of LP. 

\subsubsection{Min-Cost Network Flow for Multi-Object Tracking}
\par
\AppField{Context.}
Given detections across video frames, global data association chooses which detections belong to the same trajectories, where trajectories start and end, and which observations are false positives.

\AppField{Inputs and decision variables.}
Let $\mathcal I$ index detections ordered in time, and let $\mathcal E\subseteq\mathcal I\times\mathcal I$ contain admissible forward-time transitions.  Known costs are $\gamma_i^{\mathrm{det}}$, $\gamma_{ij}^{\mathrm{tr}}$, $\gamma_i^{\mathrm{in}}$, and $\gamma_i^{\mathrm{out}}$ for using detection $i$, transition $(i,j)$, starting at $i$, and ending at $i$.  For a specified number $M$ of tracks, let continuous relaxation variables $z_i,w_{ij},u_i,v_i\in[0,1]$ denote detection use, transition flow, source flow, and sink flow.  The min-cost flow model is
\begin{equation}
\begin{aligned}
 \min\quad &\sum_{i\in\mathcal I}\gamma_i^{\mathrm{det}}z_i
 +\sum_{(i,j)\in\mathcal E}\gamma_{ij}^{\mathrm{tr}}w_{ij} +\sum_{i\in\mathcal I}(\gamma_i^{\mathrm{in}}u_i+\gamma_i^{\mathrm{out}}v_i)\\
 \text{s.t.}\quad
 &u_i+\sum_{j:(j,i)\in\mathcal E}w_{ji}=z_i
 =v_i+\sum_{j:(i,j)\in\mathcal E}w_{ij} &&(i\in\mathcal I),\\
 &\sum_{i\in\mathcal I}u_i=M,\qquad \sum_{i\in\mathcal I}v_i=M,\\
 &0\le z_i,u_i,v_i\le1 &&(i\in\mathcal I),\\
 &0\le w_{ij}\le1 &&((i,j)\in\mathcal E).
\end{aligned}
\label{eq:tracking-flow}
\end{equation}

\AppField{Interpretation.}
The two equalities at each detection enforce flow conservation: a selected detection has exactly one incoming source/transition unit and one outgoing sink/transition unit.  Source and sink constraints create $M$ complete trajectories, and unit capacities prevent one detection from serving two trajectories.  Costs are commonly negative log-likelihood components.

\AppField{MP classification and structure.}
This is an LP with a node--arc incidence matrix.  With integral supplies and capacities, the incidence matrix is totally unimodular; therefore an extreme optimal solution is integral even though variables are declared continuous.  Additional non-network coupling, such as general pairwise exclusion or higher-order occlusion constraints, can destroy this property.

\AppField{Solution approaches.}
Successive shortest paths,  network simplex, and specialized dynamic graph algorithms exploit the structure; generic LP solvers also certify optimality.  

\AppField{MP relationship and citation.}
\textbf{Explicit formulation.}  \citet{zhang2008tracking} map maximum-a-posteriori global association to min-cost network flow and extend the network for occlusion.  Equation~\eqref{eq:tracking-flow} is a normalized core of that explicit construction.

\AppField{Limitations.}
The formulation presumes calibrated local costs, a chosen transition graph, and  a fixed track count $M$.  Long-range interactions, identity ambiguity, and errors in the detection stage are only represented through the supplied arcs and costs.

\section{Quadratic Programming}
\label{sec:qp}
QP has a quadratic objective and linear constraints.  A positive-semidefinite Hessian yields a convex problem, while indefinite curvature yields nonconvexity even without integrality.

\subsection{Local Reconstruction Weights in Locally Linear Embedding}
\AppField{Context.}
Locally Linear Embedding (LLE) estimates, for each observation, weights that reconstruct it from nearby observations.  The low-dimensional embedding is subsequently asked to preserve these local linear relations.

\AppField{Inputs and decision variables.}
For observation $i\in\mathcal I$, let $\mathbf a_i\in\R^d$ be observed, let $\mathcal N_i\subseteq\mathcal I\setminus\{i\}$ be its fixed neighborhood, and let $w_{ij}\in\R$ be the reconstruction weight on neighbor $j$.  The local program is
\begin{equation}
 \begin{aligned}
 \min_{(w_{ij})_{j\in\mathcal I}}\quad
 &\left\lVert \mathbf a_i-\sum_{j\in\mathcal I}w_{ij}\mathbf a_j\right\rVert_2^2\\
 \text{s.t.}\quad
 &\sum_{j\in\mathcal I}w_{ij}=1,\qquad
 w_{ij}=0\quad(j\notin\mathcal N_i).
 \end{aligned}
 \label{eq:lle-weights}
\end{equation}

\AppField{Interpretation.}
The objective measures local reconstruction error.  The support constraint prohibits non-neighbor influence, and the linear sum-to-one condition makes the weights invariant to translating all inputs by a common vector.

\AppField{MP classification and structure.}
Equation~\eqref{eq:lle-weights} is an equality-constrained convex QP because its  Hessian is positive semidefinite.  The problems decompose by $i$.  A singular neighborhood covariance can make weights nonunique, though the objective remains convex.

\AppField{Solution approaches.}
Each small QP can be solved through a KKT linear system; diagonal regularization stabilizes nearly singular neighborhoods.  Parallelism is immediate across observations.  Neighborhood construction, not the QP alone, often dominates large-scale cost.

\AppField{MP relationship and citation.}
\textbf{Explicit formulation.}  The constrained reconstruction problem is a defining step of LLE in \citet{roweis2000lle}; \eqref{eq:lle-weights} retains its neighborhood and linear-invariance structure.

\AppField{Limitations.}
The model treats neighborhoods as correct and fixed, is sensitive to sampling density and noise, and does not itself determine the final embedding dimension or guarantee preservation of global geometry.

\section{Binary and Mixed-Integer Programming}
\label{sec:mip}
Binary integer programming (BIP) represents yes/no and assignment choices.  Mixed-integer linear programming (MILP) adds continuous variables with linear algebra; mixed-integer quadratic programming (MIQP) permits quadratic objectives or constraints.  Branch-and-bound, cutting planes, presolve, relaxations, and primal heuristics are the main exact-solver ingredients, but problem-specific decomposition is usually decisive at ML scale.

\subsection{RAG Context Selection under a Token Budget}
\AppField{Context.}
After retrieval, a RAG system must choose a compact, nonredundant collection of passages that supports a query without exceeding the generator's context and auxiliary resource limits.

\AppField{Inputs and decision variables.}
For a fixed query $q$, let groups $\mathcal G$ partition candidate passages, let $\mathcal I_g$ be the candidates in group $g$, and let $\mathcal H$ index resources.  Candidate $i$ has known faithfulness/relevance utility $r_i$, resource consumption $c_{ih}\ge0$ for $h\in\mathcal H$ (including tokens), and each resource has budget $B_h$.  The binary decision $z_i$ indicates inclusion:
\begin{equation}
 \begin{aligned}
 \max_{z_i\in\B}\quad &\sum_{g\in\mathcal G}\sum_{i\in\mathcal I_g} r_i z_i\\
 \text{s.t.}\quad
 &\sum_{g\in\mathcal G}\sum_{i\in\mathcal I_g}c_{ih}z_i\le B_h &&(h\in\mathcal H),\\
 &\sum_{i\in\mathcal I_g}z_i\le1 &&(g\in\mathcal G).
 \end{aligned}
 \label{eq:rag-mmkp}
\end{equation}
Equality may replace the final inequality when exactly one representative per mandatory group is required.

\AppField{Interpretation.}
The objective aggregates estimated evidential value.  Resource rows enforce the token window and any redundancy or latency limits, while group-choice constraints avoid selecting near-duplicate alternatives from the same semantic group.

\AppField{MP classification and structure.}
Equation~\eqref{eq:rag-mmkp} is a multidimensional multiple-choice knapsack problem (MMKP), hence a BIP and NP-hard in general.  Its LP relaxation gives an upper bound; a single resource and no group coupling reduce it to 0--1 knapsack.

\AppField{Solution approaches.}
Dynamic programming with Pareto-state pruning, branch-and-bound, Lagrangian relaxation of budgets, and value-density heuristics are representative.  Learned utilities can be recalibrated online, but prediction uncertainty is separate from combinatorial optimality.

\AppField{MP relationship and citation.}
\textbf{Explicit formulation.}  \citet{xu2026selfcorrecting} explicitly formulate context selection as MMKP and solve it with Pareto-pruned dynamic programming before an NLI-guided tree-search generation stage.  The verifiable record is the 2026 arXiv preprint; no proceedings metadata is assumed here.

\AppField{Limitations.}
Additive utilities imperfectly represent interactions among passages, group construction may be brittle, and token counts do not capture attention dilution or ordering.  The downstream generator can still ignore selected evidence.

\subsection{Mixture-of-Experts Routing}
\AppField{Context.}
A sparse mixture-of-experts layer assigns tokens to a small number of expert subnetworks.  Routing should exploit affinity scores while preventing expert overload.

\AppField{Inputs and decision variables.}
Let $\mathcal I$ index tokens, $\mathcal K$ experts, $r_{ik}$ a known gate score, $h_i$ the required number of assignments for token $i$, and $C_k$ expert $k$'s capacity.  Binary $z_{ik}$ assigns token $i$ to expert $k$:
\begin{equation}
 \begin{aligned}
 \max_{z_{ik}\in\B}\quad &\sum_{i\in\mathcal I}\sum_{k\in\mathcal K}r_{ik}z_{ik}\\
 \text{s.t.}\quad
 &\sum_{k\in\mathcal K}z_{ik}=h_i &&(i\in\mathcal I),\\
 &\sum_{i\in\mathcal I}z_{ik}\le C_k &&(k\in\mathcal K).
 \end{aligned}
 \label{eq:moe-routing}
\end{equation}

\AppField{Interpretation.}
The first constraints provide the requested sparse computation per token; the second implement finite expert buckets.  Expert-choice routing reverses the conventional viewpoint by having experts choose high-scoring tokens, naturally fixing the second constraint while allowing heterogeneous token multiplicity.

\AppField{MP classification and structure.}
The normalized model is a capacitated bipartite BIP.  When $h_i$, $C_k$, and the constraint matrix have pure transportation structure, the LP relaxation is integral; side constraints or nonlinear gate effects may remove that property.

\AppField{Solution approaches.}
Exact min-cost flow is available for the pure assignment form.  Production routers commonly use parallel top-$k$, sorting, capacity dropping, or local balancing because solving a global assignment at every layer and batch can be too expensive.

\AppField{MP relationship and citation.}
\textbf{Natural MP reformulation.}  \citet{zhou2022expertchoice} introduce expert-choice routing with fixed expert bucket sizes and variable experts per token.  Equation~\eqref{eq:moe-routing} formalizes the capacity-aware assignment principle, but we do not claim their router invokes a BIP solver.

\AppField{Limitations.}
Scores are treated as fixed and additive, communication topology and expert specialization dynamics are omitted, and a locally optimal route need not optimize end-to-end training loss.

\subsection{Tool or Agent-Component Selection}
\AppField{Context.}
An agent composer selects which tools, models, memories, or specialized agents to activate for a task, balancing expected capability against monetary, latency, and token costs.

\AppField{Inputs and decision variables.}
Let $\mathcal I$ index components, $r_i$ be tested or predicted task utility, $c_i$ cost, and $B$ a total budget.  Let $\mathcal E^-$ be known incompatible unordered pairs and $\mathcal D$ be ordered dependency pairs, where $(i,j)\in\mathcal D$ means component $i$ requires $j$.  Binary $z_i$ selects component $i$:
\begin{equation}
 \begin{aligned}
 \max_{z_i\in\B}\quad &\sum_{i\in\mathcal I}r_i z_i\\
 \text{s.t.}\quad
 &\sum_{i\in\mathcal I}c_i z_i\le B,\\
 &z_i+z_j\le1 &&(\{i,j\}\in\mathcal E^-),\\
 &z_i\le z_j &&((i,j)\in\mathcal D).
 \end{aligned}
 \label{eq:agent-components}
\end{equation}

\AppField{Interpretation.}
The knapsack row enforces the operating budget; conflict rows prevent unusable combinations; dependency rows make a selected component bring its prerequisite.  Utilities can be updated after lightweight tests of candidates.

\AppField{MP and structure.}
Equation~\eqref{eq:agent-components} is a knapsack-style BIP with side constraints and is NP-hard.  Dependencies can sometimes be contracted; conflicts induce a packing structure.  Additive utility neglects complementarity unless interaction variables are added.

\AppField{Solution approaches.}
MIP solvers provide bounds for moderate inventories.  Online knapsack policies, greedy marginal-utility selection, Lagrangian prices, and cached evaluations support dynamic inventories and tight latency.

\AppField{MP relationship and citation.}
\textbf{Explicit formulation.}  \citet{yuan2025agentcomposition} explicitly frame automated agentic component selection through a knapsack approach that jointly accounts for utility, cost, and compatibility.  Equation~\eqref{eq:agent-components} is a compact deterministic version of that design problem.

\AppField{Limitations.}
Component utility is task- and order-dependent, tests may be noisy, and budget feasibility does not guarantee safe tool interaction.  Execution planning after selection is outside the model.

\subsection{Notification/Email Allocation and Frequency Optimization}
\AppField{Context.}
A communication platform decides which message to send to which user and when, seeking engagement or long-term value without overwhelming users or violating operational quotas.

\AppField{Inputs and decision variables.}
Let $\mathcal I$ index users, $\mathcal J$ campaigns/messages, $\mathcal T$ eligible times, and $\mathcal H$ channels.  Message $j$ uses channel $h(j)\in\mathcal H$.  Known values are predicted utility $r_{ij\tau}$, cost $c_{ij\tau}$, user cap $f_i$, channel cap $f_{ih}$, campaign lower/upper quotas $L_j,U_j$, and budget $B$.  Binary $z_{ij\tau}$ denotes a send:
\begin{equation}
 \begin{aligned}
 \max_{z_{ij\tau}\in\B}\quad &\sum_{i\in\mathcal I}\sum_{j\in\mathcal J}\sum_{\tau\in\mathcal T} r_{ij\tau}z_{ij\tau}\\
 \text{s.t.}\quad
 &\sum_{j\in\mathcal J}\sum_{\tau\in\mathcal T}z_{ij\tau}\le f_i &&(i\in\mathcal I),\\
 &\sum_{j:h(j)=h}\sum_{\tau\in\mathcal T}z_{ij\tau}\le f_{ih} &&(i\in\mathcal I,h\in\mathcal H),\\
 &L_j\le\sum_{i\in\mathcal I}\sum_{\tau\in\mathcal T}z_{ij\tau}\le U_j &&(j\in\mathcal J),\\
 &\sum_{i,j,\tau}c_{ij\tau}z_{ij\tau}\le B.
 \end{aligned}
 \label{eq:notification}
\end{equation}

\AppField{Interpretation.}
Predicted long-term value is aggregated across sends.  User and channel rows control frequency, campaign rows implement delivery commitments, and the last row limits shared resources.

\AppField{MP and structure.}
This generalized allocation is a large BIP with knapsack and quota coupling.  Without coupling, choices decompose by user; dual prices for shared constraints restore approximate decomposition.  Lower quotas can make the model infeasible and should be monitored explicitly.

\AppField{Solution approaches.}
Distributed Lagrangian methods, primal--dual pricing, LP relaxation plus repair, and partitioned heuristics are natural at billion-variable scale.  \citet{zhang2020billion} demonstrate distributed near-optimal methods for generalized industrial knapsack models, but such methods are not universally exact.

\AppField{MP relationship and citation.}
\textbf{Natural MP reformulation.}  \citet{zhao2018notification} describe a deployed ML system for per-user notification-volume control; \citet{zhang2020billion} explicitly study massive knapsack optimization in a different production setting.  Equation~\eqref{eq:notification} is the paper's joint user--message--time BIP and should not be attributed verbatim to either source.

\AppField{Limitations.}
The model treats predicted effects as additive, abstracts causal carryover and fatigue, and does not encode content safety or consent semantics.  Long-term value estimates can shift after the policy changes exposure.

\subsection{Test-Time Reasoning and Exploration-Budget Allocation}
\AppField{Context.}
An inference service allocates a finite test-time-compute pool across heterogeneous queries.  Harder queries may benefit from more samples, search depth, or verification, but allocating the largest level to every query is infeasible.

\AppField{Inputs and decision variables.}
Let $\mathcal Q$ index queries and $\mathcal K_q$ the permitted compute levels for query $q$.  Offline evaluation or a predictor supplies expected accuracy/utility $r_{qk}$ and compute cost $c_{qk}$ for level $k\in\mathcal K_q$.  Binary $z_{qk}$ chooses a level and $B$ is the batch-wide compute budget:
\begin{equation}
 \begin{aligned}
 \max_{z_{qk}\in\B}\quad &\sum_{q\in\mathcal Q}\sum_{k\in\mathcal K_q}r_{qk}z_{qk}\\
 \text{s.t.}\quad
 &\sum_{k\in\mathcal K_q}z_{qk}=1 &&(q\in\mathcal Q),\\
 &\sum_{q\in\mathcal Q}\sum_{k\in\mathcal K_q}c_{qk}z_{qk}\le B.
 \end{aligned}
 \label{eq:test-time-compute}
\end{equation}

\AppField{Interpretation.}
Each query receives exactly one inference configuration; the single coupling row makes compute fungible across the batch.  The utilities may represent pass probability, verifier score, or downstream value, but must be calibrated at each compute level.

\AppField{MP and structure.}
Equation~\eqref{eq:test-time-compute} is a multiple-choice knapsack BIP.  It is NP-hard, but admits pseudopolynomial dynamic programming for integer budgets and a strong LP relaxation when level curves are well behaved.

\AppField{Solution approaches.}
Dynamic programming, Lagrangian pricing of compute, greedy marginal-gain allocation, and batch-level MIP are viable.  Within a query, the selected level can govern tree search or repeated sampling; that execution policy is not a second application.

\AppField{MP relationship and citation.}
\textbf{Natural MP reformulation.}  \citet{snell2025testtime} show that test-time compute should be allocated adaptively according to problem difficulty and inference strategy.  They do not present the across-query BIP in \eqref{eq:test-time-compute}; it is our resource-allocation abstraction of their compute-optimal principle.

\AppField{Limitations.}
Expected gains can be nonstationary and correlated across samples, compute levels may be interruptible rather than discrete, and verifier error can make measured utility diverge from answer correctness.

\subsection{Optimal Decision-Tree Construction}
\AppField{Context.}
Optimal-tree learning selects all splits and leaf predictions jointly, trading training errors against tree complexity rather than greedily committing to one split at a time.

\AppField{Inputs and decision variables.}
Let $\mathcal I$ index labeled observations $(\mathbf a_i,y_i)$ with $\mathbf a_i\in\R^d$ and class $y_i\in\mathcal C$.  A fixed tree template has branch nodes $\mathcal B$ and candidate leaves $\mathcal L$.  For leaf $l$, sets $\mathcal A_l^L,\mathcal A_l^R\subseteq\mathcal B$ identify left and right ancestors.  Binary $d_t$ activates split $t$, $p_{jt}$ selects feature $j$, $z_{il}$ assigns sample $i$ to leaf $l$, $u_l$ activates leaf $l$, $q_{lc}$ predicts class $c$, and $e_i$ records error; $b_t\in\R$ is a threshold.  With valid feature bounds $U_t$, big-$M$ value $M_t$, strict-separation constant $\varepsilon>0$, and complexity price $\lambda\ge0$, a compact MILP is
\begin{align}
 \min\quad &\sum_{i\in\mathcal I}e_i+\lambda\sum_{t\in\mathcal B}d_t
 \label{eq:optimal-tree}\\
 \text{s.t.}\quad
 &\sum_{j=1}^d p_{jt}=d_t,\quad 0\le b_t\le U_t d_t &&(t\in\mathcal B),\nonumber\\
 &\sum_{l\in\mathcal L}z_{il}=1 &&(i\in\mathcal I),\nonumber\\
 &\sum_{c\in\mathcal C}q_{lc}=u_l,\quad z_{il}\le u_l &&(i\in\mathcal I,l\in\mathcal L),\nonumber\\
 &z_{il}\le d_t &&(i\in\mathcal I,l\in\mathcal L,t\in\mathcal A_l^L\cup\mathcal A_l^R),\nonumber\\
 &\sum_{j=1}^d a_{ij}p_{jt}\le b_t+M_t(1-z_{il})
 &&(i,l,t\in\mathcal A_l^L),\nonumber\\
 &\sum_{j=1}^d a_{ij}p_{jt}\ge b_t+\varepsilon-M_t(1-z_{il})
 &&(i,l,t\in\mathcal A_l^R),\nonumber\\
 &e_i\ge z_{il}-q_{l,y_i} &&(i\in\mathcal I,l\in\mathcal L),\nonumber\\
 &d_t,p_{jt},z_{il},u_l,q_{lc},e_i\in\B.\nonumber
\end{align}
Here and below, a quantifier such as $(i,l,t\in\mathcal A_l^L)$ abbreviates all $i\in\mathcal I$, $l\in\mathcal L$, and $t\in\mathcal A_l^L$.

\AppField{Interpretation.}
Feature-selection and threshold rows define univariate splits.  Assignment and path rows route every sample consistently; leaf rows choose one prediction at each used leaf; the last inequalities detect misclassification.  Complexity regularization favors fewer active splits.  Standard hierarchy constraints, omitted only to save display width, prevent an inactive parent from having an active descendant.

\AppField{MP and structure.}
The model is an MILP.  Integrality represents the tree topology and routing, while big-$M$ path constraints can weaken the relaxation.  The problem is combinatorial and globally optimizing even shallow trees can be expensive.

\AppField{Solution approaches.}
Branch-and-cut with feature/threshold preprocessing, warm starts from greedy trees, symmetry breaking, valid inequalities, and decomposition by leaves are representative.  Early termination yields a feasible tree plus an optimality gap rather than a false claim of exactness.

\AppField{MP relationship and citation.}
\textbf{Explicit formulation.}  \citet{bertsimas2017oct} formulate optimal classification trees using mixed-integer optimization and jointly optimize accuracy and complexity.  Equation~\eqref{eq:optimal-tree} is a normalized compact template, not a line-for-line reproduction.

\AppField{Limitations.}
Big-$M$ constants require careful scaling, training accuracy may not predict calibration or robustness, and exact solution is generally restricted to moderate data or shallow templates.

\subsection{Best Subset Selection}
\AppField{Context.}
Best subset regression chooses at most $m$ features while fitting coefficients by squared error, producing a sparse predictive model with an explicit cardinality budget.

\AppField{Inputs and decision variables.}
Let $\mathbf A\in\R^{n\times d}$ be the design matrix and $\mathbf y\in\R^n$ the response.  Continuous $\theta_j\in\R$ are regression coefficients, binary $z_j$ selects feature $j$, and valid bounds $M_j>0$ link the two:
\begin{equation}
 \begin{aligned}
 \min_{\theta\in\R^d,\,z\in\B^d}\quad
 &\frac12\lVert\mathbf y-\mathbf A\theta\rVert_2^2\\
 \text{s.t.}\quad
 &-M_jz_j\le\theta_j\le M_jz_j &&(j=1,\ldots,d),\\
 &\sum_{j=1}^d z_j\le m.
 \end{aligned}
 \label{eq:best-subset}
\end{equation}

\AppField{Interpretation.}
If $z_j=0$, the linking bounds force coefficient $\theta_j$ to zero; if $z_j=1$, it may vary within a valid range.  The final row enforces the feature budget and the objective measures residual fit.

\AppField{MP and structure.}
Equation~\eqref{eq:best-subset} is a convex MIQP: its continuous quadratic objective has a positive-semidefinite Hessian, but binary selection makes the full problem nonconvex.  Loose $M_j$ values degrade the relaxation; perspective or indicator formulations can strengthen it.

\AppField{Solution approaches.}
Modern MIP solvers combine branch-and-bound, cutting planes, and warm starts from sparse heuristics.  Discrete first-order methods, screening, and local swaps rapidly produce incumbents; only a solver bound certifies their global gap.

\AppField{MP relationship and citation.}
\textbf{Explicit formulation.}  \citet{bertsimas2016subset} develop a modern mixed-integer optimization treatment of best subset selection, including warm starts and certifiable solutions.

\AppField{Limitations.}
The formulation assumes squared-error adequacy and valid coefficient bounds, can be unstable under collinearity, and does not automatically account for post-selection uncertainty.

\subsection{Constrained Maximum-Inner-Product Retrieval}
\AppField{Context.}
Constrained maximum-inner-product search (MIPS) retrieves a top-$m$ set for a query while enforcing category, fairness, or coverage requirements.  Pairwise diversity is an extension of the same retrieval application, not a separate application.

\AppField{Inputs and decision variables.}
For query $q$, let $\mathbf v_q\in\R^d$ be its embedding and $\mathbf a_i\in\R^d$ item $i$'s embedding; the known score is $r_i=\langle\mathbf v_q,\mathbf a_i\rangle$.  Matrix $\mathbf C\in\R^{h\times|\mathcal I|}$ and vectors $\mathbf l,\mathbf u\in\R^h$ encode $h$ lower/upper category, fairness, or coverage bounds.  Binary $z_i$ selects an item.  The constrained linear case and pairwise-diverse case are respectively
\begin{align}
 \max_{z\in\B^{|\mathcal I|}}\quad &\sum_{i\in\mathcal I}r_i z_i
 &&\text{s.t.}\quad \sum_i z_i=m,\quad \mathbf l\le\mathbf C z\le\mathbf u,
 \label{eq:constrained-mips}\\
 \max_{z\in\B^{|\mathcal I|}}\quad
 &\sum_{i\in\mathcal I}r_i z_i-\lambda\sum_{i<j}s_{ij}z_i z_j
 &&\text{s.t.}\quad \sum_i z_i=m,\quad \mathbf l\le\mathbf C z\le\mathbf u,
 \label{eq:diverse-mips}
\end{align}
where $s_{ij}\ge0$ is pairwise similarity and $\lambda\ge0$ its penalty weight.

\AppField{Interpretation.}
Equation~\eqref{eq:constrained-mips} chooses the highest total inner-product score among feasible lists.  Equation~\eqref{eq:diverse-mips} subtracts redundancy: two similar items incur a penalty only when both are selected.

\AppField{MP and structure.}
The linear model is a BIP/MILP.  The pairwise term makes the extension a binary quadratic program (BQP), commonly handled as MIQP after retaining the product or as MILP after exact product linearization.  With nonnegative $s_{ij}$ and maximization of its negative, the continuous quadratic curvature need not make a standard convex MIQP representation; discreteness remains central.

\AppField{Solution approaches.}
Candidate generation by approximate MIPS followed by exact constrained reranking is common.  Branch-and-bound, Lagrangian relaxation of group bounds, local swaps, and diversity-aware approximate indexes trade candidate recall against optimization cost.

\AppField{MP relationship and citation.}
\textbf{Natural MP reformulation.}  \citet{hirata2023categorical} define categorical diversity-aware inner-product search and give a sublinear approximation with a probabilistic success guarantee.  Equations~\eqref{eq:constrained-mips}--\eqref{eq:diverse-mips} are a general MP normalization, particularly for linear fairness/coverage and pairwise diversity; they are not attributed as the paper's solver model.

\AppField{Limitations.}
The model optimizes only over the candidate universe supplied to it, similarity can be a poor proxy for semantic redundancy, and list constraints may conflict.  Approximate indexing can miss the globally optimal feasible set before optimization starts.

\subsection{Mixed-Precision Neural-Network Quantization}
\AppField{Context.}
Mixed-precision quantization assigns a bit-width to each network layer, minimizing predicted accuracy degradation under memory or latency limits.

\AppField{Inputs and decision variables.}
Let $\mathcal I$ index layers and $\mathcal K_i$ the allowed bit-width options for layer $i$.  Known $d_{ik}\ge0$ is a sensitivity-based loss proxy and $c_{ikh}\ge0$ is consumption of resource $h\in\mathcal H$, whose budget is $B_h$.  Binary $z_{ik}$ selects option $k$:
\begin{equation}
 \begin{aligned}
 \min_{z_{ik}\in\B}\quad &\sum_{i\in\mathcal I}\sum_{k\in\mathcal K_i}d_{ik}z_{ik}\\
 \text{s.t.}\quad
 &\sum_{k\in\mathcal K_i}z_{ik}=1 &&(i\in\mathcal I),\\
 &\sum_{i\in\mathcal I}\sum_{k\in\mathcal K_i}c_{ikh}z_{ik}\le B_h &&(h\in\mathcal H).
 \end{aligned}
 \label{eq:mixed-precision}
\end{equation}

\AppField{Interpretation.}
Exactly one precision is chosen per layer.  Sensitivity scores approximate model perturbation, and resource rows enforce size, energy, or measured latency limits.

\AppField{MP and structure.}
This is a multiple-choice multidimensional BIP/ILP and is NP-hard.  It has strong separability except for the few budget rows, making Lagrangian relaxation attractive.

\AppField{Solution approaches.}
ILP solvers, dynamic programming for a single discretized budget, dual pricing, and sensitivity-guided rounding are representative.  Hardware measurements should replace nominal bit counts when latency is the actual constraint.

\AppField{MP relationship and citation.}
\textbf{Explicit formulation.}  \citet{yao2021hawq} explicitly solve an integer linear program that balances Hessian-based perturbation against hardware-aware memory and latency constraints.

\AppField{Limitations.}
Layerwise sensitivities neglect cross-layer interactions, latency is device- and kernel-specific, and accuracy after fine-tuning may not follow the static proxy.

\subsection{Diversity-Aware Recommendation}
\AppField{Context.}
Recommendation-list construction chooses a fixed-size set whose members are individually relevant yet collectively nonmonotonous.

\AppField{Inputs and decision variables.}
For a user, candidate item $i\in\mathcal I$ has relevance $r_i$, pairwise similarity $s_{ij}\ge0$, and binary selection $z_i$.  Optional linear feasibility data $\mathbf C,\mathbf u$ are known and locally defined.  With trade-off $\lambda\ge0$,
\begin{equation}
 \begin{aligned}
 \max_{z\in\B^{|\mathcal I|}}\quad
 &\sum_{i\in\mathcal I}r_i z_i-\lambda\sum_{i<j}s_{ij}z_i z_j\\
 \text{s.t.}\quad
 &\sum_{i\in\mathcal I}z_i=m,\qquad \mathbf C z\le\mathbf u.
 \end{aligned}
 \label{eq:diverse-recommendation}
\end{equation}

\AppField{Interpretation.}
The first term rewards predicted preference; the second penalizes similar co-selected pairs.  Cardinality fixes list length, and optional linear rows can represent inventory or policy feasibility.

\AppField{MP and structure.}
Equation~\eqref{eq:diverse-recommendation} is a BQP/MIQP-style model.  Binary products can be linearized with auxiliary variables and McCormick inequalities, producing an MILP at quadratic size.  Depending on the similarity matrix and sign convention, the continuous relaxation may be nonconcave for maximization.

\AppField{Solution approaches.}
Exact MIQP/MILP is possible after candidate pruning; greedy marginal-gain selection, local search, and maximal-marginal-relevance-style reranking scale farther but do not automatically certify optimality.

\AppField{MP relationship and citation.}
\textbf{Natural MP reformulation.}  \citet{zhang2008monotony} formalize the accuracy--diversity trade-off for recommendation lists and propose diversification methods.  Equation~\eqref{eq:diverse-recommendation} is a canonical binary quadratic rendering of that trade-off, not a claim about the exact solver used in the paper.

\AppField{Limitations.}
Pairwise similarity is only one notion of diversity, the scalar $\lambda$ hides user-specific trade-offs, and relevance predictions are affected by exposure and position not modeled here.

\subsection{Globally Optimized $k$-Means or Minimum Sum-of-Squares Clustering}
\AppField{Context.}
Minimum sum-of-squares clustering jointly assigns observations to $K$ clusters and estimates centroids, seeking a global solution rather than a local fixed point of alternating updates.

\AppField{Inputs and decision variables.}
Let $\mathbf a_i\in\R^d$ be observed for $i\in\mathcal I$, and let $\mathcal K=\{1,\ldots,K\}$.  Binary $z_{ik}$ assigns observation $i$ to cluster $k$, continuous $\bm\mu_k\in\R^d$ is its centroid, and $t_{ik}\in\R_+$ is an assignment cost.  Suppose all feasible centroids lie in known bounds that yield valid constants $M_{ik}$.  An exact mixed-integer quadratic-constraint model is
\begin{equation}
 \begin{aligned}
 \min_{z,\bm\mu,t}\quad &\sum_{i\in\mathcal I}\sum_{k\in\mathcal K}t_{ik}\\
 \text{s.t.}\quad
 &\sum_{k\in\mathcal K}z_{ik}=1 &&(i\in\mathcal I),\\
 &\sum_{i\in\mathcal I}z_{ik}\ge1 &&(k\in\mathcal K),\\
 &t_{ik}\ge\lVert\mathbf a_i-\bm\mu_k\rVert_2^2-M_{ik}(1-z_{ik}) &&(i\in\mathcal I,k\in\mathcal K),\\
 &z_{ik}\in\B,\quad t_{ik}\ge0,\quad \bm\mu_k\in\R^d.
 \end{aligned}
 \label{eq:global-kmeans}
\end{equation}

\AppField{Interpretation.}
Every observation joins exactly one nonempty cluster.  If $z_{ik}=1$, the epigraph constraint charges its squared distance to centroid $k$; a valid $M_{ik}$ deactivates the constraint otherwise.

\AppField{MP and structure.}
Directly writing $\sum_{ik}z_{ik}\lVert\mathbf a_i-\bm\mu_k\rVert_2^2$ gives a nonconvex mixed-integer nonlinear program (MINLP), not a BQP.  Equation~\eqref{eq:global-kmeans} is a mixed-integer convex quadratic-constraint reformulation when bounds are valid.  Other exact developments produce MIQP or MILP models through disjunctions, distance variables, and linearization.  Cluster-label symmetry and weak big-$M$ bounds are major obstacles.

\AppField{Solution approaches.}
Spatial branch-and-bound, branch-and-cut, symmetry breaking, bound tightening, semidefinite or assignment relaxations, and incumbents from repeated local clustering are representative.  Tailored propagation and cuts improve general-purpose MINLP solvers, but global certification remains difficult.

\AppField{MP relationship and citation.}
\textbf{Explicit formulation.}  \citet{burgard2023clustering} develop MIP techniques for the minimum sum-of-squares problem, while \citet{agoston2024kmeans} give MILP formulations.  These sources support the distinctions among MINLP, quadratic, and reformulated linear models.

\AppField{Limitations.}
The squared Euclidean objective favors spherical clusters, $K$ is fixed, outliers can dominate, and global methods remain limited to far smaller instances than local heuristics.

\subsection{Video Frame or Keyframe Selection}
\AppField{Context.}
Long-video systems select a small set of query-relevant frames that covers distinct events and fits a multimodal model's visual-token budget.

\AppField{Inputs and decision variables.}
Let $\mathcal I$ index frames in temporal order, $r_i$ be query relevance, $s_{ij}\ge0$ visual/semantic redundancy, $m$ the frame limit, and $\mathcal E_\Delta$ pairs that are too close in time to co-select.  Binary $z_i$ selects a frame and $\lambda\ge0$ controls redundancy:
\begin{equation}
 \begin{aligned}
 \max_{z\in\B^{|\mathcal I|}}\quad
 &\sum_{i\in\mathcal I}r_i z_i-\lambda\sum_{i<j}s_{ij}z_i z_j\\
 \text{s.t.}\quad
 &\sum_{i\in\mathcal I}z_i\le m,\\
 &z_i+z_j\le1 &&((i,j)\in\mathcal E_\Delta).
 \end{aligned}
 \label{eq:keyframe}
\end{equation}

\AppField{Interpretation.}
The objective balances query relevance and content diversity.  The first row limits multimodal compute; optional temporal-separation rows prevent near-duplicate adjacent frames.

\AppField{MP and structure.}
This is an integer quadratic program/BQP.  Pairwise terms make exact optimization NP-hard; an MILP linearization uses one auxiliary product variable per retained pair.  Sparsifying the similarity graph reduces both memory and solve time.

\AppField{Solution approaches.}
MIQP/MILP handles short candidate lists; greedy marginal-gain search, local swaps, temporal segmentation, and learned candidate pruning support long videos.  A heuristic should be reported as such even when the objective is exact.

\AppField{MP relationship and citation.}
\textbf{Explicit formulation.}  \citet{fang2026keyframe} explicitly formulate query-relevant, diversity-aware keyframe choice as integer quadratic programming and use a customized greedy alternative for efficiency.

\AppField{Limitations.}
Framewise relevance can miss actions requiring temporal context, pairwise redundancy does not ensure narrative continuity, and downstream generated narratives can introduce errors not represented by \eqref{eq:keyframe}.

\section{Conic Programming}
\label{sec:conic}
Conic programs optimize a linear objective over linear constraints and membership in convex cones.  Second-order-cone programs (SOCPs) capture norm epigraphs; semidefinite programs (SDPs) optimize over the positive-semidefinite cone.  Conic duality supplies strong bounds, but generic SDP scaling is substantially less favorable than LP or SOCP scaling.

\subsection{Group Lasso as SOCP}
\AppField{Context.}
Group Lasso selects or removes predefined blocks of coefficients together, matching settings in which predictors form meaningful factors.

\AppField{Inputs and decision variables.}
Let $\mathbf A\in\R^{n\times d}$ and $\mathbf y\in\R^n$ be observed, and let the feature partition be $\{\mathcal J_g:g\in\mathcal G\}$.  Coefficients in group $g$ form $\bm\beta_g\in\R^{|\mathcal J_g|}$, with known weight $\omega_g>0$ and penalty $\lambda\ge0$.  The Group Lasso objective is
\begin{equation}
 \min_{\bm\beta\in\R^d}\quad
 \frac12\lVert\mathbf y-\mathbf A\bm\beta\rVert_2^2
 +\lambda\sum_{g\in\mathcal G}\omega_g\lVert\bm\beta_g\rVert_2.
 \label{eq:group-lasso}
\end{equation}
Define the rotated cone $\mathcal Q_r^3=\{(a,b,c)\in\R_+\times\R_+\times\R:2ab\ge c^2\}$.  With epigraph variables $u,v,t_g\ge0$, an equivalent SOCP is
\begin{equation}
 \begin{aligned}
 \min_{\bm\beta,u,v,t}\quad &v+\lambda\sum_{g\in\mathcal G}\omega_g t_g\\
 \text{s.t.}\quad
 &\lVert\mathbf y-\mathbf A\bm\beta\rVert_2\le u,\\
 &\lVert\bm\beta_g\rVert_2\le t_g &&(g\in\mathcal G),\\
 &(v,1,u)\in\mathcal Q_r^3.
 \end{aligned}
 \label{eq:group-lasso-socp}
\end{equation}

\AppField{Interpretation.}
The first cone bounds residual norm, group cones bound coefficient norms, and the rotated cone enforces $v\ge u^2/2$.  At optimum all relevant epigraph bounds are tight, recovering \eqref{eq:group-lasso}.

\AppField{MP and structure.}
Equation~\eqref{eq:group-lasso-socp} is a convex SOCP.  Nondifferentiability at a zero group produces block sparsity without integrality.  Strict convexity, and therefore coefficient uniqueness, depends on the design and loss.

\AppField{Solution approaches.}
Interior-point SOCP solvers are effective at moderate scale.  Block coordinate descent, proximal gradient, and screening exploit separability for larger designs; they solve the original composite objective without explicitly forming cones.

\AppField{MP relationship and citation.}
\textbf{Equivalent formulation.}  \citet{yuan2006grouplasso} introduce grouped variable selection through the Group Lasso.  The explicit rotated-cone representation in \eqref{eq:group-lasso-socp} is mathematically equivalent but is not attributed as their implementation.

\AppField{Limitations.}
Groups must be specified, the penalty can over-shrink large or correlated groups, and the basic model does not capture overlapping hierarchy unless additional variables are introduced.

\subsection{Sparse PCA as SDP}
\AppField{Context.}
Sparse principal component analysis (PCA) seeks a direction explaining high variance while using few coordinates, improving interpretability at the cost of a combinatorial sparsity constraint.

\AppField{Inputs and decision variables.}
Let $\bm\Sigma\in\R^{d\times d}$ be a known sample covariance matrix with $\bm\Sigma\succeq0$.  A loading $\mathbf v\in\R^d$ induces the lifted decision $\mathbf X=\mathbf v\mathbf v^\top$.  For sparsity budget $\kappa>0$, the rank-one model and its convex relaxation are
\begin{align}
 \max_{\mathbf X}\quad &\langle\bm\Sigma,\mathbf X\rangle
 &&\text{s.t.}\quad \tr(\mathbf X)=1, \lVert\mathbf X\rVert_1\le\kappa,
 \ \mathbf X\succeq0, \operatorname{rank}(\mathbf X)=1,
 \label{eq:sparse-pca-rank1}\\
 \max_{\mathbf X}\quad &\langle\bm\Sigma,\mathbf X\rangle
 &&\text{s.t.}\quad \tr(\mathbf X)=1, \lVert\mathbf X\rVert_1\le\kappa,
 \ \mathbf X\succeq0,
 \label{eq:sparse-pca-sdp}
\end{align}
where $\lVert\mathbf X\rVert_1=\sum_{i=1}^d\sum_{j=1}^d|X_{ij}|$ is the entrywise $\ell_1$ norm.

\AppField{Interpretation.}
Trace normalization corresponds to $\lVert\mathbf v\rVert_2=1$ in the rank-one model; explained variance is $\mathbf v^\top\bm\Sigma\mathbf v=\langle\bm\Sigma,\mathbf X\rangle$.  The entrywise budget promotes a loading supported on few coordinates.  Dropping rank one enlarges the feasible set and gives an upper bound.

\AppField{MP and structure.}
Equation~\eqref{eq:sparse-pca-rank1} is nonconvex because of rank.  Equation~\eqref{eq:sparse-pca-sdp} is an SDP after linearizing absolute values.  Its solution need not be rank one, so a loading recovered from its leading eigenvector is generally approximate.

\AppField{Solution approaches.}
Interior-point SDP works at small to medium dimension; first-order spectral methods, smoothing, low-rank factorization, and screening scale farther.  Randomized or leading-eigenvector rounding converts a higher-rank relaxation into a sparse loading but may lose the SDP bound's value.

\AppField{MP relationship and citation.}
\textbf{Explicit formulation.}  \citet{daspremont2004sparsepca} directly derive an SDP relaxation for sparse PCA from a rank-constrained lifted formulation.

\AppField{Limitations.}
The $d\times d$ lift has quadratic storage, relaxation solutions can have high rank, and a covariance estimate contaminated by noise can dominate any optimization guarantee.

\section{Bilevel Programming}
\label{sec:bilevel}
Bilevel programming nests one optimization problem inside another.  With outer feasible set $\Lambda$, inner parameter set $\Theta$, upper objective $\Phi$, and lower objective $\Psi$, the normalized form is
\begin{equation}
 \min_{\lambda\in\Lambda}\ \Phi\bigl(\lambda,\theta^\star(\lambda)\bigr)
 \quad\text{s.t.}\quad
 \theta^\star(\lambda)\in\arg\min_{\theta\in\Theta}\Psi(\theta,\lambda).
 \label{eq:bilevel-common}
\end{equation}
When the inner minimizer is not unique, optimistic and pessimistic conventions differ.  Hypergradients can be obtained by unrolling, implicit differentiation, or approximation, but nonconvex inner learning means these procedures need not solve the global bilevel problem.

\subsection{Hyperparameter Optimization}
\AppField{Context.}
Hyperparameter optimization selects regularization, data-processing, or training-control parameters for validation performance after the model has been trained with those choices.

\AppField{Inputs and decision variables.}
Training and validation samples define losses $L_{\mathrm{tr}}(\theta,\lambda)$ and $L_{\mathrm{val}}(\theta,\lambda)$.  The outer decision is $\lambda\in\Lambda$ and the learned inner decision is $\theta\in\Theta$:
\begin{equation}
 \min_{\lambda\in\Lambda}\ L_{\mathrm{val}}\bigl(\theta^\star(\lambda),\lambda\bigr)
 \quad\text{s.t.}\quad
 \theta^\star(\lambda)\in\arg\min_{\theta\in\Theta}L_{\mathrm{tr}}(\theta,\lambda).
 \label{eq:hpo}
\end{equation}

\AppField{Interpretation.}
The inner problem trains a model for the proposed hyperparameters; the outer objective measures generalization on held-out data.  Separating the data prevents the outer level from merely choosing settings that overfit the training loss.

\AppField{MP and structure.}
This is bilevel programming.  It is smooth only under strong differentiability and solution-regularity conditions; deep-network instances are nonconvex at both practical and theoretical levels.

\AppField{Solution approaches.}
Unrolled differentiation backpropagates through a finite training trajectory; implicit differentiation solves a linear system involving the inner Hessian; truncated or first-order hypergradients reduce memory.  Black-box search is useful for discrete or nonsmooth $\lambda$ but does not exploit the bilevel derivative.

\AppField{MP relationship and citation.}
\textbf{Explicit formulation.}  \citet{franceschi2018bilevel} explicitly unify gradient-based hyperparameter optimization and meta-learning through bilevel programming and analyze finite inner approximations.

\AppField{Limitations.}
Validation reuse can overfit, inner training is rarely solved exactly, and hypergradient quality depends on optimizer dynamics, stopping rules, and Hessian approximations.

\subsection{Neural Architecture Search}
\AppField{Context.}
Neural architecture search (NAS) chooses operations and connectivity using validation performance, while network weights are learned on training data.

\AppField{Inputs and decision variables.}
Let $\lambda$ collect architecture weights and let $\Lambda$ be a product of simplices: for every edge $e$ in a fixed search graph and operation $k\in\mathcal K_e$, $\lambda_{ek}\ge0$ and $\sum_{k\in\mathcal K_e}\lambda_{ek}=1$.  Let $\theta$ be network weights.  With architecture-dependent losses,
\begin{equation}
 \begin{aligned}
 \min_{\lambda\in\Lambda}\quad &L_{\mathrm{val}}\bigl(\theta^\star(\lambda),\lambda\bigr)\\
 \text{s.t.}\quad
 &\theta^\star(\lambda)\in\arg\min_{\theta\in\Theta}L_{\mathrm{tr}}(\theta,\lambda),\\
 &\Lambda=\left\{\lambda:\lambda_{ek}\ge0,\ \sum_{k\in\mathcal K_e}\lambda_{ek}=1\ \text{for every }e\right\}.
 \end{aligned}
 \label{eq:nas}
\end{equation}

\AppField{Interpretation.}
Soft architecture weights mix candidate operations so gradients can update architecture and model parameters.  A discrete architecture is derived after search, typically by retaining the largest weight on each edge.

\AppField{MP and structure.}
Equation~\eqref{eq:nas} is a nonconvex bilevel program with a continuous relaxation of an underlying discrete design.  The relaxed optimum and the discretized architecture need not agree.

\AppField{Solution approaches.}
Alternating inner/outer gradient steps, one-step unrolling, implicit gradients, and weight sharing reduce search cost.  Retraining the selected discrete architecture is necessary to evaluate the final design independently of shared weights.

\AppField{MP relationship and citation.}
\textbf{Explicit formulation.}  \citet{liu2019darts} explicitly formulate differentiable architecture search as a bilevel problem and introduce the continuous operation relaxation represented in \eqref{eq:nas}.

\AppField{Limitations.}
Weight sharing biases architecture comparisons, short inner optimization changes the intended objective, discretization creates a gap, and search results can be sensitive to random initialization and regularization.

\subsection{Coreset Selection and Data Subsampling}
\AppField{Context.}
Coreset selection chooses a small weighted or probabilistic training subset whose trained model performs well on the full distribution or a validation set.

\AppField{Inputs and decision variables.}
Let training samples $(\mathbf a_i,y_i)$ be indexed by $\mathcal I$.  Outer decisions $w_i\in[0,1]$ are selection probabilities (binary indicators are an exact-subset alternative), with expected size at most $m$.  For per-sample loss $\ell_i(\theta)=\ell(\theta;(\mathbf a_i,y_i))$, define
\begin{equation}
 \begin{aligned}
 \min_{w\in[0,1]^{|\mathcal I|}}\quad
 &L_{\mathrm{val}}\bigl(\theta^\star(w)\bigr)\\
 \text{s.t.}\quad
 &\theta^\star(w)\in\arg\min_{\theta\in\Theta}
 \frac{1}{\sum_i w_i}\sum_{i\in\mathcal I}w_i\ell_i(\theta),\\
 &\sum_{i\in\mathcal I}w_i\le m,\qquad \sum_i w_i>0.
 \end{aligned}
 \label{eq:coreset-bilevel}
\end{equation}

\AppField{Interpretation.}
Selection weights define the data seen by inner training; outer validation evaluates the resulting model.  The size row limits expected coreset mass.  Probabilistic sampling turns a discrete outer search into stochastic continuous optimization.

\AppField{MP and structure.}
The model is bilevel and generally nonconvex.  Binary selection is combinatorial and nondifferentiable; probabilistic weights permit unbiased score-function gradients but can have high variance.

\AppField{Solution approaches.}
Policy-gradient estimators, implicit or unrolled hypergradients for continuous weights, warm-started inner training, and class-balanced sampling are representative.  Rounding or sampling produces the final subset.

\AppField{MP relationship and citation.}
\textbf{Explicit formulation.}  \citet{zhou2022coreset} explicitly propose probabilistic bilevel coreset selection and an unbiased policy-gradient solver, including the size-controlled sampling interpretation.

\AppField{Limitations.}
Repeated inner training is expensive, validation choice controls what the coreset preserves, probabilities do not directly guarantee one realized subset's quality, and nonconvex learning complicates reproducibility.

\subsection{Autonomous Scientific Discovery and AI Scientist Systems}
\AppField{Context.}
An autonomous discovery system may optimize not only candidate hypotheses or designs, but also the scientific objective or computable evaluation criteria used to judge them.

\AppField{Inputs and decision variables.}
Let outer decision $\lambda\in\Lambda$ encode a computable scientific criterion or its weights, let candidate design $\mathbf x\in\mathcal X$ include discrete hypotheses and continuous physical parameters, let $\Psi(\mathbf x,\lambda)$ be the criterion-driven search score, and let $\Phi(\lambda,\mathbf x)$ be an external scientific evaluation obtained from simulation, experiment, or a broader goal.  A normative bilevel abstraction is
\begin{equation}
 \min_{\lambda\in\Lambda}\ \Phi\bigl(\lambda,\mathbf x^\star(\lambda)\bigr)
 \quad\text{s.t.}\quad
 \mathbf x^\star(\lambda)\in\arg\min_{\mathbf x\in\mathcal X}\Psi(\mathbf x,\lambda).
 \label{eq:scientific-bilevel}
\end{equation}

\AppField{Interpretation.}
The inner loop searches for a candidate under the current scientific goal; the outer loop analyzes outcomes and revises the goal or scoring function.  In simulation-grounded variants, language models may propose discrete structures while differentiable simulation tunes continuous parameters.

\AppField{MP and structure.}
Equation~\eqref{eq:scientific-bilevel} is a bilevel representation, but its components may be discrete, black-box, stochastic, and nonstationary.  An agentic outer loop that writes new objectives is not automatically equivalent to conventional differentiable bilevel programming.

\AppField{Solution approaches.}
Iterative propose--simulate--analyze loops, surrogate modeling, derivative-based inner simulation, discrete search, and human checkpoints can be combined.  Solver convergence claims apply only to the mathematical subproblem actually solved, not to open-ended scientific validity.

\AppField{MP relationship and citation.}
\textbf{Equivalent formulation.}  \citet{ma2024scientific} explicitly call their LLM--simulation Scientific Generative Agent a bilevel optimization framework.  \citet{du2025saga} describe SAGA's outer goal-evolving agents and inner optimizer as a bi-level architecture.  Equation~\eqref{eq:scientific-bilevel} faithfully normalizes the architecture, while avoiding a claim that SAGA is a smooth or globally solved bilevel program.

\AppField{Limitations.}
Scientific objectives are incomplete proxies, simulators can be misspecified, experiments can be costly, and autonomous criterion revision creates auditability and safety challenges beyond optimization.

\subsection{Diffusion-Model Hyperparameter and Noise-Schedule Optimization}
\AppField{Context.}
Diffusion systems choose a noise schedule or fine-tuning hyperparameters whose quality can be assessed only after training or running an inner generative process.

\AppField{Inputs and decision variables.}
Let $\lambda\in\Lambda$ parameterize a feasible noise schedule, let $\theta$ be denoiser parameters, let $\tau\in\{1,\ldots,T\}$ be a sampled diffusion time, and let $\bm\epsilon\sim P_\epsilon$ be noise.  A clean datum is $\mathbf a\sim P_{\mathrm{data}}$, the schedule-generated noisy state is $\mathbf x_\tau(\mathbf a,\bm\epsilon;\lambda)$, and the denoiser is $\bm\epsilon_\theta(\cdot,\tau)$.  For an outer generation-quality loss $L_{\mathrm{gen}}$, the bilevel model is
\begin{equation}
 \begin{aligned}
 \min_{\lambda\in\Lambda}\quad
 &L_{\mathrm{gen}}\bigl(\theta^\star(\lambda),\lambda\bigr)\\
 \text{s.t.}\quad
 &\theta^\star(\lambda)\in\arg\min_{\theta\in\Theta}
 \mathbb E_{\substack{\mathbf a\sim P_{\mathrm{data}},\,\tau,\\\bm\epsilon\sim P_\epsilon}}
 \left[\left\lVert\bm\epsilon-\bm\epsilon_\theta(\mathbf x_\tau(\mathbf a,\bm\epsilon;\lambda),\tau)\right\rVert_2^2\right].
 \end{aligned}
 \label{eq:diffusion-bilevel}
\end{equation}

\AppField{Interpretation.}
The inner level trains the denoiser for the proposed schedule; the outer level scores generated validation samples or a differentiable quality surrogate.  In an inference-only fine-tuning regime, the inner map can instead be a sampling or adaptation solver.

\AppField{MP and structure.}
This is stochastic generative bilevel optimization.  The inner problem is nonconvex and its output is a distribution induced by iterative sampling, so standard finite-dimensional implicit differentiation can be computationally prohibitive.

\AppField{Solution approaches.}
First-order hypergradient estimators, reparameterization of training dynamics, sample-efficient outer estimators, truncated unrolling, and schedule parameterization reduce cost.  Each method optimizes an approximation determined by finite samples and finite inner work.

\AppField{MP relationship and citation.}
\textbf{Explicit formulation.}  \citet{xiao2025diffusion} explicitly formalize diffusion fine-tuning and noise-schedule training as generative bilevel optimization and derive first-order estimators for both regimes.

\AppField{Limitations.}
Generation metrics are imperfect and costly, stochastic gradient variance is high, and a schedule optimized for one data distribution, sampler, or compute budget may not transfer.

\section{Multi-Objective Optimization}
\label{sec:multiobjective}
Multi-objective optimization retains several objectives rather than collapsing them prematurely to one scalar.  Dominance and Pareto stationarity replace a unique total order, while scalarization or common-descent subproblems provide computational steps.

\subsection{Multi-Task Learning}
\AppField{Context.}
Multi-task learning shares representations across tasks whose losses can conflict.  The decision is a common parameterization that is not improvable on one task without degrading at least one other task.

\AppField{Inputs and decision variables.}
Let $\mathcal T=\{1,\ldots,T\}$ index tasks, let $L_t(\theta)$ be differentiable task $t$'s loss, and let $\theta\in\Theta$ contain shared and task-specific model parameters.  The vector problem is
\begin{equation}
 \min_{\theta\in\Theta}\quad
 \bigl(L_1(\theta),\ldots,L_T(\theta)\bigr).
 \label{eq:mtl-vector}
\end{equation}
At a current $\theta$, let $\mathbf g_t=\nabla_\theta L_t(\theta)\in\R^p$.  Simplex weights $\alpha_t$ solve the common-descent QP
\begin{equation}
 \begin{aligned}
 \min_{\alpha\in\R^T}\quad
 &\frac12\left\lVert\sum_{t=1}^{T}\alpha_t\mathbf g_t\right\rVert_2^2\\
 \text{s.t.}\quad &\alpha_t\ge0\ (t=1,\ldots,T),\qquad \sum_{t=1}^{T}\alpha_t=1.
 \end{aligned}
 \label{eq:mtl-simplex-qp}
\end{equation}

\AppField{Interpretation.}
A parameter is Pareto optimal if no feasible parameter weakly improves every $L_t$ and strictly improves one.  The minimum-norm convex combination of task gradients gives a common descent direction when nonzero; zero indicates first-order Pareto stationarity under standard regularity.

\AppField{MP and structure.}
Equation~\eqref{eq:mtl-vector} is multi-objective and generally nonconvex for neural models.  Equation~\eqref{eq:mtl-simplex-qp} is a small convex QP in $T$ variables even when $p$, the parameter dimension, is large.

\AppField{Solution approaches.}
Gradient-based multiple-gradient descent, weighted scalarizations, preference-conditioned training \citet{momma2022multi}, and approximate low-dimensional gradient representations are representative.  Solving one QP produces one Pareto-stationary trajectory, not the entire Pareto front.

\AppField{MP relationship and citation.}
\textbf{Explicit formulation.}  \citet{sener2018multitask} explicitly cast multi-task learning as multi-objective optimization and use a simplex-constrained minimum-norm subproblem, with a scalable upper-bound construction for deep networks.

\AppField{Limitations.}
Pareto stationarity is not global Pareto optimality in nonconvex models, gradient scales influence trade-offs, and a technically Pareto solution may still be unacceptable on a critical task.

\section{Inverse Optimization}
\label{sec:inverse}
Inverse optimization observes decisions and infers objective parameters under which those decisions are optimal or nearly optimal.  Identifiability is central: many objectives can rationalize the same behavior.

\subsection{Inverse Reinforcement Learning}
\AppField{Context.}
Inverse reinforcement learning (IRL) infers a reward from expert behavior rather than learning a policy from a supplied reward.

\AppField{Inputs and decision variables.}
Consider a known finite discounted Markov decision process with states $\mathcal S$, actions $\mathcal A$, transition probabilities $P(s'\mid s,a)$, discount $\gamma\in[0,1)$, and demonstrated deterministic policy $\pi_E(s)$.  Known features $\bm\phi(s,a)\in\R^d$ define reward $r_\theta(s,a)=\theta^\top\bm\phi(s,a)$.  Decisions are reward weights $\theta\in\R^d$, expert-policy values $V_s\in\R$, and slacks $\xi_{sa}\ge0$.  Define the continuation operator $(\mathcal P_aV)_s=\sum_{s'\in\mathcal S}P(s'\mid s,a)V_{s'}$.  For regularization $C>0$, a maximum-margin inverse-optimality QP is
\begin{equation}
\begin{aligned}
 \min_{\theta,V,\xi}\quad &\frac12\lVert\theta\rVert_2^2+C\sum_{s\in\mathcal S}\sum_{a\ne\pi_E(s)}\xi_{sa}
 \\
 \text{s.t.}\quad
 &V_s=r_\theta(s,\pi_E(s))+\gamma(\mathcal P_{\pi_E(s)}V)_s &&(s\in\mathcal S),\\
 &V_s\ge r_\theta(s,a)+\gamma(\mathcal P_aV)_s+1-\xi_{sa}
 &&(s\in\mathcal S,a\ne\pi_E(s)),\\
 &\xi_{sa}\ge0 &&(s\in\mathcal S,a\ne\pi_E(s)).
\end{aligned}
\label{eq:irl}
\end{equation}

\AppField{Interpretation.}
The first row evaluates the demonstrated policy.  The second asks its action to beat each alternative by unit margin, softened by $\xi_{sa}$.  Regularization selects one explanation among otherwise scale-ambiguous rewards.

\AppField{MP and structure.}
Equation~\eqref{eq:irl} is a convex QP with linear Bellman inverse-optimality constraints.  Alternative IRL formulations use feature-expectation separation or a bilevel condition that places the expert policy in an inner control argmax.

\AppField{Solution approaches.}
QP/LP solvers apply for known finite dynamics.  Constraint generation avoids enumerating all policies, while forward control oracles find violated optimality conditions.  Approximate dynamics and large state spaces require sampling or function approximation.

\AppField{MP relationship and citation.}
\textbf{Explicit formulation.}  \citet{ng2000irl, dong2018generalized, dong2024towards} formulate IRL through linear inequalities expressing policy optimality and develop algorithms that resolve reward degeneracy through margins and priors.  Equation~\eqref{eq:irl} is a normalized soft-margin Bellman version of that explicit inverse-optimization idea.

\AppField{Limitations.}
Rewards are generally nonidentifiable up to transformations, demonstrations may be suboptimal, dynamics may be unknown, and inferred intent need not correspond to ethically acceptable goals.

\section{Distributionally Robust Optimization}
\label{sec:dro}
DRO protects against probability laws near an empirical distribution.  For parameter $\theta\in\Theta$, random observation $\xi$, loss $\ell(\theta;\xi)$, empirical law $\widehat P$, and ambiguity set $\mathcal U(\widehat P)$, the common structure is
\begin{equation}
 \min_{\theta\in\Theta}\ \max_{Q\in\mathcal U(\widehat P)}
 \mathbb E_{\xi\sim Q}[\ell(\theta;\xi)].
 \label{eq:dro-common}
\end{equation}
The ambiguity geometry determines which shifts are protected: reweighting under a divergence ball differs from transporting probability mass under a Wasserstein ball.  Duality often replaces the inner distribution problem by a finite-dimensional regularized loss.

\subsection{Contrastive Learning as DRO}
\AppField{Context.}
Contrastive learning pulls an anchor toward a positive view while separating it from negatives.  Negative sampling can be biased, so robustness is taken over plausible negative distributions.

\AppField{Inputs and decision variables.}
Let $(\mathbf a,\mathbf a^+)$ be an anchor--positive pair drawn from known empirical law $\widehat P_+$, let $\widehat P_-$ be the empirical negative law, and let encoder $f_\theta$ have parameters $\theta\in\Theta$.  Define similarity $s_\theta(\mathbf u,\mathbf v)=\langle f_\theta(\mathbf u),f_\theta(\mathbf v)\rangle$, temperature $\eta>0$, divergence $D_\varphi$, and radius $\rho\ge0$.  With
$\mathcal U_\varphi(\widehat P_-)=\{Q:D_\varphi(Q\Vert\widehat P_-)\le\rho\}$, a robust pairwise contrastive model is
\begin{equation}
 \min_{\theta\in\Theta}\ 
 \mathbb E_{(\mathbf a,\mathbf a^+)\sim\widehat P_+}
 \left[
 \max_{Q\in\mathcal U_\varphi(\widehat P_-)}
 \mathbb E_{\mathbf a^-\sim Q}
 \log\!\left(1+\exp\!\frac{s_\theta(\mathbf a,\mathbf a^-)-s_\theta(\mathbf a,\mathbf a^+)}{\eta}\right)
 \right].
 \label{eq:contrastive-dro}
\end{equation}

\AppField{Interpretation.}
The inner distribution upweights negatives that most challenge the current representation while remaining in a divergence neighborhood of observed negatives.  Radius $\rho$ controls conservatism, and temperature $\eta$ has a dual role related to robustness regularization.

\AppField{MP and structure.}
Equation~\eqref{eq:contrastive-dro} is DRO/min--max optimization.  Convex duality may solve the distributional inner problem for fixed losses, but neural dependence on $\theta$ makes the outer problem nonconvex.  Overly large ambiguity sets can emphasize false negatives or outliers.

\AppField{Solution approaches.}
Divergence dualization, adversarial reweighting of minibatch negatives, stochastic gradient descent--ascent, and adjusted contrastive losses are representative.  Finite negative batches approximate both the empirical law and its adversary.

\AppField{MP relationship and citation.}
\textbf{Explicit formulation.}  \citet{wu2023contrastive} explicitly interpret contrastive learning as DRO over the negative-sampling distribution and relate temperature to ambiguity-set control.

\AppField{Limitations.}
The ambiguity set may not match semantic shift, adversarial emphasis can amplify mislabeled negatives, and finite-batch dual estimates can be biased or unstable.

\subsection{Distributionally Robust DPO for LLM Alignment}
\AppField{Context.}
Direct preference optimization (DPO) aligns an LLM from preferred/rejected response pairs.  DRO protects against demographic, linguistic, geographic, or temporal shifts in the preference distribution.

\AppField{Inputs and decision variables.}
An observation is $\xi=(\mathbf a,y^+,y^-)$, containing prompt $\mathbf a$, preferred response $y^+$, and rejected response $y^-$.  Let $\widehat P$ be their empirical law, $\pi_\theta$ the learned policy, $\pi_{\mathrm{ref}}$ a fixed reference policy, $\beta>0$ a scale, and $\sigma(t)=(1+e^{-t})^{-1}$.  Define
\begin{equation}
 \ell_{\mathrm{DPO}}(\theta;\xi)=
 -\log\sigma\!\left(\beta\left[
 \log\frac{\pi_\theta(y^+\mid\mathbf a)}{\pi_{\mathrm{ref}}(y^+\mid\mathbf a)}-
 \log\frac{\pi_\theta(y^-\mid\mathbf a)}{\pi_{\mathrm{ref}}(y^-\mid\mathbf a)}
 \right]\right).
 \label{eq:dpo-loss}
\end{equation}
For radius $\rho\ge0$ and either divergence/transport discrepancy $D$, let $\mathcal U(\widehat P)=\{Q:D(Q,\widehat P)\le\rho\}$.  Robust DPO is
\begin{equation}
 \min_{\theta\in\Theta}\ \max_{Q\in\mathcal U(\widehat P)}
 \mathbb E_{\xi\sim Q}\bigl[\ell_{\mathrm{DPO}}(\theta;\xi)\bigr].
 \label{eq:robust-dpo}
\end{equation}

\AppField{Interpretation.}
The inner distribution concentrates on preference examples or nearby transported examples with high DPO loss.  The outer policy must improve preferred-response odds relative to the reference under this worst plausible law.

\AppField{MP and structure.}
Equation~\eqref{eq:robust-dpo} is DRO with a nonconvex neural outer problem.  Kullback--Leibler (KL) ambiguity produces reweighting; Wasserstein ambiguity permits geometric movement according to a chosen ground cost.  The ambiguity set, not only its radius, defines the protected shift.

\AppField{Solution approaches.}
Dual approximations convert the inner problem into a regularized loss, followed by stochastic gradients.  Alternating adversarial weights and policy steps, clipping, and minibatch approximations support scale but solve approximate minimax objectives.

\AppField{MP relationship and citation.}
\textbf{Explicit formulation.}  \citet{xu2025robustdpo} explicitly develop Wasserstein DPO and KL-DPO, characterize sample complexity, and derive scalable approximations for their minimax losses.

\AppField{Limitations.}
Preference labels can be systematically wrong, the transport cost may not reflect human similarity, conservatism can reduce average utility, and robustness to represented shifts is not universal alignment.

\section{Submodular Optimization}
\label{sec:submodular}
For a finite ground set $\mathcal V$, the common model is
\begin{equation}
 \max_{S\subseteq\mathcal V}\ F(S)
 \quad\text{s.t.}\quad S\in\mathcal F,
 \label{eq:submodular-common}
\end{equation}
where $\mathcal F$ is a feasible family.  A set function is normalized if $F(\varnothing)=0$, monotone if $F(A)\le F(B)$ for $A\subseteq B$, and submodular if
$F(A\cup\{i\})-F(A)\ge F(B\cup\{i\})-F(B)$ for all $A\subseteq B\subseteq\mathcal V$ and $i\notin B$.  The last property is diminishing returns.  Greedy selection attains $1-1/e$ of optimum for normalized monotone submodular maximization under a cardinality constraint; other constraints require different guarantees.

\subsection{Data Subset Selection and Active Learning}
\AppField{Context.}
A learner selects a representative training subset or, after uncertainty filtering, a batch whose labels should be acquired.  These are two modes of one approved data-selection application.

\AppField{Inputs and decision variables.}
Let $\mathcal V$ index data, $s_{ij}\ge0$ be similarity from item $i$ to representative $j$, and $m$ be the selection limit.  For $S\subseteq\mathcal V$, define facility-location value $F(S)=\sum_{i\in\mathcal V}\max_{j\in S}s_{ij}$ with $\max_{j\in\varnothing}s_{ij}=0$.  The model is
\begin{equation}
 \max_{S\subseteq\mathcal V}\quad \sum_{i\in\mathcal V}\max_{j\in S}s_{ij}
 \qquad\text{s.t.}\qquad |S|\le m.
 \label{eq:data-submodular}
\end{equation}

\AppField{Interpretation.}
Every data point is credited according to its most similar selected representative, so additional representatives have diminishing coverage benefit.  In active learning, $\mathcal V$ may first be restricted to uncertain examples, after which \eqref{eq:data-submodular} removes redundancy.

\AppField{MP and structure.}
The facility-location function is normalized, monotone, and submodular for nonnegative similarities.  Cardinality-constrained maximization is NP-hard, but standard greedy has a $1-1/e$ guarantee under exactly these assumptions.

\AppField{Solution approaches.}
Lazy greedy, stochastic greedy, memoized marginal gains, and distributed partitioning scale selection.  The guarantee concerns the chosen set function, not downstream test accuracy.

\AppField{MP relationship and citation.}
\textbf{Explicit formulation.}  \citet{wei2015submodularity} explicitly connect classifier data likelihoods to constrained submodular maximization and combine uncertainty filtering with submodular selection for active learning.

\AppField{Limitations.}
Similarity determines what ``representative'' means, uncertainty estimates may be miscalibrated, and submodular coverage does not preserve rare behavior unless encoded in the ground set or objective.

\subsection{Document Summarization}
\AppField{Context.}
Extractive summarization selects sentences that cover salient concepts and span distinct content while respecting a word budget.

\AppField{Inputs and decision variables.}
Let $\mathcal V$ index candidate sentences, $c_i$ be sentence length, $B$ the word budget, and $\mathcal G$ concepts with weights $\omega_g\ge0$ and coverage $a_{ig}\ge0$.  Let topic groups $\mathcal V_k$ partition sentences, let $r_i\ge0$ measure salience, and let $\lambda\ge0$.  The selected set $S$ solves
\begin{equation}
 \begin{aligned}
 \max_{S\subseteq\mathcal V}\quad
 &F(S)=\sum_{g\in\mathcal G}\omega_g\min\!\left\{1,\sum_{i\in S}a_{ig}\right\}
 +\lambda\sum_k\sqrt{\sum_{i\in S\cap\mathcal V_k}r_i}\\
 \text{s.t.}\quad &\sum_{i\in S}c_i\le B.
 \end{aligned}
 \label{eq:summarization}
\end{equation}

\AppField{Interpretation.}
The capped coverage term rewards mentioning important concepts without repeatedly crediting them.  The concave-over-modular topic term rewards distributing selected salience across topics.  Sentence lengths enforce the summary budget.

\AppField{MP and structure.}
With nonnegative data, $F$ is normalized, monotone, and submodular.  The length constraint is knapsack rather than pure cardinality; modified greedy/enumeration algorithms give constant-factor guarantees under appropriate assumptions.

\AppField{Solution approaches.}
Lazy cost-benefit greedy, partial enumeration plus greedy, and streaming approximations are representative.  Marginal gains can be cached by concept and topic.

\AppField{MP relationship and citation.}
\textbf{Explicit formulation.}  \citet{lin2011summarization} design monotone submodular coverage/diversity functions for document summarization and use budgeted maximization with approximation guarantees.

\AppField{Limitations.}
Extracted sentences can be incoherent, concept weights may miss factual dependencies, word cost is only a length proxy, and redundancy beyond the chosen features can remain.

\subsection{KV-Cache Eviction for LLM Inference}
\AppField{Context.}
During autoregressive LLM generation, the key--value (KV) cache grows with sequence length.  An eviction policy must retain a fixed number of tokens using only information available so far.

\AppField{Inputs and decision variables.}
At decoding step $\tau$, let $S_{\tau-1}$ be the retained token indices and $G_\tau=S_{\tau-1}\cup\{\tau\}$ the candidates after the new token arrives.  Let $o_{\tau i}\ge0$ be token $i$'s accumulated attention-derived heavy-hitter score, let $h:\R_+\to\R_+$ be nondecreasing and concave, and define $F_\tau(S)=h(\sum_{i\in S}o_{\tau i})$.  Let $R_\tau\subseteq G_\tau$ be a reserved recent window and $B$ the cache budget.  The one-step dynamic set problem is
\begin{equation}
 \begin{aligned}
 S_\tau\in\arg\max_{S\subseteq G_\tau}\quad &F_\tau(S)\\
 \text{s.t.}\quad &|S|=B,\qquad |S\setminus S_{\tau-1}|\le1,\qquad R_\tau\subseteq S.
 \end{aligned}
 \label{eq:h2o}
\end{equation}
When the cache is not full, $|S|=B$ is replaced by $|S|\le B$ and the new token is added.

\AppField{Interpretation.}
The score favors historical heavy hitters, the reserve protects recent tokens before enough attention evidence accumulates, and the dynamic constraint allows at most the arriving token to enter.  Operationally, one evaluates which single eviction from $G_\tau$ leaves the greatest score; this is not ordinary static top-$k$ because scores and the feasible ground set evolve with generation.

\AppField{MP and structure.}
Concave-over-modular $F_\tau$ is monotone submodular.  The cited method develops a broader dynamic submodular framework and a greedy eviction guarantee under its assumptions.  With the linear choice $h(t)=t$, the one-step score is modular, while the dynamic analysis and score updates remain essential.

\AppField{Solution approaches.}
H2O accumulates observed attention, reserves recent entries, and performs low-cost greedy eviction.  Exact enumeration of one removal is cheap per step; the systems challenge is obtaining and updating attention scores without offsetting memory savings.

\AppField{MP relationship and citation.}
\textbf{Explicit formulation.}  \citet{zhang2023h2o} explicitly formulate KV-cache eviction as a dynamic submodular problem, motivate heavy hitters, and retain a balance of heavy and recent tokens.  Equation~\eqref{eq:h2o} makes that evolving-set structure explicit rather than reducing the method to static top-$k$.

\AppField{Limitations.}
Past attention may not predict future need, heads and layers can disagree about importance, score collection has overhead, and evicted information cannot be recovered without recomputation.

\section{Min--Max and Saddle-Point Optimization}
\label{sec:minmax}
Min--max problems model strategic or adversarial interaction.  Convex--concave games admit strong saddle theory, but neural games are usually nonconvex--nonconcave; alternating gradients can cycle or converge to points without equilibrium guarantees.

\subsection{Generative Adversarial Networks}
\AppField{Context.}
A generative adversarial network (GAN) trains a generator to imitate data while a discriminator learns to distinguish generated from real samples.

\AppField{Inputs and decision variables.}
Let $P_{\mathrm{data}}$ be the real-data law, $P_\epsilon$ a noise law, generator $G_{\theta_G}$ have parameters $\theta_G$, and discriminator $D_{\theta_D}:\R^d\to(0,1)$ have parameters $\theta_D$.  For real sample $\mathbf a$ and noise $\bm\epsilon$, the canonical game is
\begin{equation}
 \min_{\theta_G}\max_{\theta_D}
 \left[
 \mathbb E_{\mathbf a\sim P_{\mathrm{data}}}\log D_{\theta_D}(\mathbf a)
 +\mathbb E_{\bm\epsilon\sim P_\epsilon}
 \log\bigl(1-D_{\theta_D}(G_{\theta_G}(\bm\epsilon))\bigr)
 \right].
 \label{eq:gan}
\end{equation}

\AppField{Interpretation.}
The discriminator maximizes correct log-likelihood on real and generated samples.  The generator minimizes the same value by producing samples that the discriminator accepts.  In the ideal unrestricted game, equilibrium occurs when the generated law matches the data law and the discriminator outputs one half.

\AppField{MP and structure.}
Equation~\eqref{eq:gan} is a min--max/saddle-point problem.  With neural parameterizations it is nonconvex in $\theta_G$ and nonconcave in $\theta_D$, so classical minimax interchange and global convergence conditions do not hold.

\AppField{Solution approaches.}
Alternating stochastic gradients, multiple discriminator steps, extra-gradient updates, regularization, and alternative generator losses address numerical behavior.  These are practical game-optimization methods, not exact saddle solvers.

\AppField{MP relationship and citation.}
\textbf{Explicit formulation.}  \citet{goodfellow2014gan} introduce the adversarial framework explicitly as the two-player minimax game in \eqref{eq:gan} and analyze its ideal equilibrium.

\AppField{Limitations.}
Finite models and samples invalidate the ideal analysis, training can cycle or collapse modes, and discriminator loss may correlate poorly with perceptual quality.

\section{Cross-Paradigm Comparison and Practical Considerations}
\label{sec:comparison}
Table~\ref{tab:comparison} compares modeling consequences rather than declaring one paradigm universally preferable.  A system can cross rows: prediction may be trained by bilevel or DRO methods and then consumed by an MILP; an SDP can bound a discrete model; a Lagrangian dual can decompose a BIP into per-request decisions.

\begin{table*}
\caption{Cross-paradigm comparison.  Complexity statements are typical worst-case characterizations, not instance-level runtime promises.}
\label{tab:comparison}
\scriptsize
\setlength{\tabcolsep}{2pt}
\begin{tabularx}{\textwidth}{P{1.35cm}P{1.35cm}P{1.45cm}P{1.45cm}P{1.9cm}P{1.9cm}Y}
\toprule
Paradigm & Variables & Convexity & Typical complexity & Common solution methods & Scalability bottleneck & Representative ML/AI role\\
\midrule
LP & Continuous; sometimes integral by structure & Convex & Polynomial-time algorithms; often very fast in practice & Simplex, interior point, network algorithms, decomposition & Rows/columns and communication & Sparse recovery; global association\\
QP & Continuous & Convex for PSD Hessian; otherwise nonconvex & Polynomial for convex QP; hard in general & KKT solves, active set, first order & Dense Hessians and conditioning & Local reconstruction weights\\
BIP/MIP & Binary plus optional continuous & Nonconvex feasible set & NP-hard in general & Branch-and-cut, relaxation, decomposition, heuristics & Search tree and weak bounds & Selection, assignment, routing, structured models\\
Conic & Continuous matrices/vectors & Convex for standard cones & Polynomial-time theory & Interior point, first-order splitting, spectral methods & SDP matrix storage/factorization & Group norms; lifted sparse components\\
Bilevel & Outer and inner, continuous or discrete & Usually nonconvex & Hard even for simple subclasses & Unrolling, implicit differentiation, nested search & Repeated inner solves and hypergradient memory & Model/configuration design\\
Multi-\newline objective & Usually continuous parameters & Inherits objective geometry & Pareto set may be large & Scalarization, common-descent QP, continuation & Many tasks/objectives and gradient storage & Conflicting task losses\\
Inverse & Objective parameters plus forward variables & Convex only for selected models & Depends on forward problem and identifiability & KKT/Bellman constraints, constraint generation & Repeated forward-oracle calls & Recover latent preferences/rewards\\
DRO & Model plus adversarial distribution & Convex in selected loss/ambiguity pairs; neural outer nonconvex & Inner problem can be infinite-dimensional before duality & Dual reformulation, reweighting, stochastic descent--ascent & Adversarial sampling and ambiguity calibration & Robust representations/alignment\\
Submodular & Set-valued/\newline discrete & Discrete diminishing returns & NP-hard maximization; approximable under structure & Greedy, lazy/stochastic greedy, local search & Marginal-gain evaluation and ground-set size & Coverage/diversity under budgets\\
Min--max & Two or more players & Favorable if convex--concave; neural games are not & Global neural equilibria generally intractable & Alternating gradients, extragradient, regularization & Instability and unequal player progress & Adversarial generation\\
\bottomrule
\end{tabularx}
\end{table*}

\subsection{From Predictions to Coefficients}
Scores such as $r_i$, sensitivities, costs, and similarities are usually estimated.  Optimization then treats them as data, so a perfect solver can still make a poor decision when coefficients are biased or shifted.  Calibration, uncertainty sets, causal evaluation, and post-decision monitoring are therefore part of the optimization pipeline.  The relationship label is equally important: a natural MP reformulation can be valuable for counterfactual analysis even if a latency-critical implementation uses sorting or greedy search.

\subsection{Relaxations, Certificates, and Approximation}
An LP relaxation of a BIP supplies a bound, not necessarily a deployable decision.  Total unimodularity in \eqref{eq:tracking-flow} is exceptional because it closes the integrality gap under stated conditions.  SDP relaxation in \eqref{eq:sparse-pca-sdp} similarly gives a bound, but recovering a rank-one loading can lose optimality.  Submodular greedy guarantees compare against the optimum of the specified set function; they do not validate that function as a proxy for human utility.  For nonconvex bilevel and min--max learning, a stationary point is weaker than a global certificate and should be described accordingly.

\subsection{Decomposition and Scale}
Many large models are separable except for a few budgets: notification decisions couple through quotas, quantization choices through hardware limits, and compute levels through a batch budget.  Lagrangian multipliers turn shared constraints into prices and expose parallel subproblems.  Column generation is useful when feasible actions are numerous but only a few enter the master problem.  Candidate generation before constrained reranking is indispensable in retrieval, but its recall ceiling must be distinguished from the optimizer's gap on the reduced candidate set.

\subsection{Choosing a Formulation}
Modeling should begin with the decision and guarantee required.  Use network LP when incidence structure is genuine; do not add binaries merely because the underlying objects are discrete.  Use MIP when logical feasibility and certificates matter at manageable scale.  Use SOCP/SDP when norm or matrix geometry provides a useful convex representation or bound.  Use bilevel models when one learned response is conditional on another decision, DRO when the uncertainty is distributional, submodularity when coverage exhibits diminishing returns, and min--max models when there are genuinely opposed players.  Hybrid algorithms are normal, but the reported guarantee must correspond to the part actually solved.

\section{Research Challenges and Opportunities}
\label{sec:challenges}
\paragraph{End-to-end uncertainty.}
Most formulations separate prediction from optimization.  A central challenge is to propagate epistemic and distributional uncertainty in relevance, cost, and response models into decisions without making robust models unusably conservative.

\paragraph{Differentiation through discrete and approximate solvers.}
Learning coefficients end to end requires gradients through LP/MIP layers, greedy subset selection, or truncated inner training.  Useful estimators must reconcile discontinuity, bias, memory, and solver tolerances while preserving feasibility.

\paragraph{Scalable certificates.}
Information-access instances can contain billions of candidate actions.  Better decompositions, learned primal heuristics, safe screening, and distributed bounds could make certificates available beyond small offline models.  A learned heuristic should accelerate search without being allowed to invalidate exact pruning.

\paragraph{Dynamic objectives and constraints.}
LLM inference, notification systems, and autonomous discovery operate in nonstationary environments.  Static optimality is inadequate when actions alter future data, preferences, or scientific goals.  Online and dynamic formulations need regret or stability statements aligned with the actual feedback loop.

\paragraph{Semantics of diversity and fairness.}
Pairwise similarity, category quotas, and coverage features are convenient proxies.  Research should connect them to user-level outcomes, intersectional constraints, and causal exposure effects while retaining algorithms that scale to retrieval and recommendation workloads.

\paragraph{Auditability of agentic optimization.}
When agents revise objectives or compose tools, reproducibility requires versioned objectives, constraint provenance, solver status, and traces of approximations.  This is a modeling opportunity: provenance and safety conditions can become first-class constraints rather than post hoc documentation.

\paragraph{Benchmarking formulation fidelity.}
Comparisons should separate coefficient error, candidate-generation loss, relaxation gap, optimization gap, and execution drift.  Reporting only final task quality obscures whether the mathematical model, its learned inputs, or its solver is responsible for failure.


\bibliographystyle{ACM-Reference-Format}
\bibliography{references}

\end{document}